\documentclass{article}

\usepackage[utf8]{inputenc}
\usepackage[T1]{fontenc}  
\usepackage{lmodern}
\usepackage[english]{babel}

\usepackage[a4paper]{geometry}
\usepackage{amsmath,amsthm,amssymb}
\usepackage{bm}

\usepackage{hyperref}
\usepackage{csquotes}

\usepackage{authblk}
\usepackage[title]{appendix}
\usepackage{siunitx}

\usepackage[backend=biber,firstinits=true,
maxbibnames=15,sorting=none]{biblatex}
\usepackage{tikz,pgfplots,pgfplotstable}
\pgfplotsset{compat=newest}
\usetikzlibrary{arrows,snakes,backgrounds,positioning,calc}

\usepackage{subfig}
\def\letters{a,b,c,d,e,f,g,h,i,j,k,l,m,n,o,p,q,r,s,t,u,v,w,x,y,z}
\def\Letters{A,B,C,D,E,F,G,H,I,J,K,L,M,N,O,P,Q,R,S,T,U,V,W,X,Y,Z}
\makeatletter
\@for \@l:=\Letters \do{%
  \expandafter\edef\csname\@l bb\endcsname{\noexpand\ensuremath{\noexpand\mathbb{\@l}}}%
  \expandafter\edef\csname\@l bf\endcsname{{\noexpand\bf \@l}}%
  \expandafter\edef\csname\@l cal\endcsname{\noexpand\ensuremath{\noexpand\mathcal{\@l}}}%
  \expandafter\edef\csname\@l eu\endcsname{\noexpand\ensuremath{\noexpand\EuScript{\@l}}}%
  \expandafter\edef\csname\@l frak\endcsname{\noexpand\ensuremath{\noexpand\mathfrak{\@l}}}%
  \expandafter\edef\csname\@l rm\endcsname{{\noexpand\rm \@l}}%
  \expandafter\edef\csname\@l scr\endcsname{\noexpand\ensuremath{\noexpand\mathscr{\@l}}}%
}
\@for \@l:=\letters \do{%
  \expandafter\edef\csname\@l bf\endcsname{{\noexpand\bf \@l}}%
  \expandafter\edef\csname\@l frak\endcsname{\noexpand\ensuremath{\noexpand\mathfrak{\@l}}}%
  \expandafter\edef\csname\@l scr\endcsname{\noexpand\ensuremath{\noexpand\mathscr{\@l}}}%
}
\makeatother

\newcommand{\chibf}{\boldsymbol{\chi}}
\newcommand{\Sigmaref}{\Sigma_\mathrm{ref}}
\newcommand{\Qref}{Q_\mathrm{ref}}

\newcommand{\isdef}{\mathrel{\mathrel{\mathop:}=}}
\newcommand{\defis}{\mathrel{=\mathrel{\mathop:}}}

\newcommand{\dd}{\operatorname{d}\!}

\newcommand{\Cov}{\operatorname{Cov}}
\DeclareMathOperator{\StdDev}{StdDev}

\newcommand{\xmat}{\widehat{\xbf}}
\newcommand{\zmat}{\widehat{\zbf}}
\newcommand{\rparam}{{\boldsymbol\xi}}
\newcommand{\rparamscal}{{\xi}}

\DeclareMathOperator*{\argmin}{\operatorname{argmin}}
\newcommand{\refd}{{\operatorname*{ref}}}

\title{Space-time shape uncertainties in the forward and inverse problem of electrocardiography}

\author[]{Lia Gander}
\author[]{Rolf Krause}
\author[]{Michael Multerer}
\author[]{Simone Pezzuto\thanks{Corresponding Author: \url{simone.pezzuto@usi.ch}}}
\affil[]{Center for Computational Medicine in Cardiology,\par
Euler Institute, \par
Universit\`a della Svizzera italiana,\par
via la Santa 1, CH-6900 Lugano, Switzerland}

\begin{document}

\maketitle

\begin{abstract}

In electrocardiography, the ``classic'' inverse problem is the reconstruction
of electric potentials at a surface enclosing the heart from remote recordings
at the body surface and an accurate description of the anatomy. The latter being
affected by noise and obtained with limited resolution due to clinical constraints,
a possibly large uncertainty may be perpetuated in the inverse reconstruction.

The purpose of this work is to study the effect of shape uncertainty on the
forward and the inverse problem of electrocardiography.  To this aim, the
problem is first recast into a boundary integral formulation and then
discretised with a collocation method to achieve high convergence rates and a
fast time to solution. The shape uncertainty of the domain is represented by a
random deformation field defined on a reference configuration. We propose a
periodic-in-time covariance kernel for the random field and approximate the
Karhunen-Lo\`eve expansion using low-rank techniques for fast sampling. The
space-time uncertainty in the expected potential and its variance is evaluated
with an anisotropic sparse quadrature approach and validated by a quasi-Monte
Carlo method.

We present several numerical experiments on a simplified but physiologically
grounded 2-dimensional geometry to illustrate the validity of the
approach. The tested parametric dimension ranged from 100 up to 600.
For the forward problem the sparse quadrature
is very effective.  In the inverse problem, the sparse quadrature and the
quasi-Monte Carlo method perform as expected, except for the total variation
regularisation, where convergence is limited by lack of regularity.  We
finally investigate an $H^{1/2}$ regularisation, which naturally stems
from the boundary integral formulation, and compare it to more classical
approaches.

\medskip\noindent
\textit{Keywords}: space-time shape uncertainty; boundary integral formulation;
sparse quadrature; quasi-Monte Carlo method; inverse problem
of electrocardiography; $H^{1/2}$ regularisation

\end{abstract}

\section{Introduction}

Electrocardiographic recordings at the body surface, such as the standard
12-lead electrocardiogram (ECG), are a direct consequence of the electric
activity of the heart.  The spatial resolution of such recordings depends
on the number of electrodes placed on the chest, ranging from no more
than ten of the standard ECG to hundreds of electrodes composing
high-density body surface potential maps (BSPMs). Along with an accurate
description of the torso anatomy, BSPMs are sufficiently informative to
enable a non-invasive characterization of cardiac electrophysiology,
termed ECG imaging~\cite{Ramanathan2004}. ECG imaging technologies have
been extensively validated in experimental, animal and, more recently,
clinical settings, with promising results~\cite{Cluitmans2018}.

The ECG imaging problem can mathematically be formulated as an inverse problem.
The most classical formulation of it is associated with a potential-based
forward problem, which amounts to determining the BSPM on the chest from the
electric potential at a surface enclosing the heart, e.g., the epicardium.
The inverse problem of electrocardiography consists therefore in recovering the
epicardial potential from the BSPM~\cite{colli1985inverse}. As typical for
inverse problems, however, the ECG imaging problem is severely ill-posed in the
sense of Hadamard, since arbitrarily small perturbations of the BSPM, such as
noise, may yield large variations in the epicardial potential. It thus
requires regularisation to penalize non-physical solutions and a strategy
to optimally select the associated regularisation parameter. For the ECG inverse
problem, several kinds of regularisations have been proposed over the last
three decades, see \cite{karoui2018evaluation} for a review.  The standard
Tikhonov regularisations of zeroth, first or second
order~\cite{colli1985inverse,rudy1988inverse} are easy to implement thanks to the
closed-form solution of the quadratic inverse problem. Closely related to the
Tikhonov regularisation is the generalized truncated singular value
decomposition, in which small singular values of the forward operator are
filtered out~\cite{hansen1992modified}. The inverse problem can also be
interpreted from a Bayesian perspective~\cite{serinagaoglu2005bayesian}. In
this case the maximum a posteriori estimator, which matches the classical
Tikhonov solution under the usual hypothesis of a standard Gaussian prior
distribution, offers more flexibility in embodying prior knowledge in the
inverse problem, e.g., from a training data-set. Non-smooth regularisations
such as Total Variation (TV) are a valid alternative to quadratic approaches
in the presence of sharp gradients in the reconstructed
data~\cite{ghosh2009application}, but lead to a significantly more difficult
solution of the inverse problem. Approximated or smoothed versions of TV are
therefore popular~\cite{karoui2018evaluation}. More recently, physiology-based
and spatio-temporal regularisation approaches have also been
considered~\cite{cluitmans2017physiology,Schuler2018}.

The potential-based forward problem is particularly attractive when the torso
is assumed as a homogeneous electric conductor. Then, the forward problem
can be conveniently recast into an integral formulation involving only the
boundaries of the torso, that is the epicardium and the chest,
with no need of solving the problem in the full 3-dimensional domain~\cite{barr1977relating,kress1989linear}. The numerical
treatment of the boundary formulation is however less practical than a
standard finite element approach in 3-D, as it requires special care in the
treatment of singular boundary integrals on piecewise smooth surfaces, like
triangulated surfaces obtained from the segmentation of cardiac images.  
Acquired with given clinical constraints, imaging
data are typically of limited resolution and noisy. Therefore, the
segmentation of selected cavities is challenging and certainly subject to
uncertainty. In the case of the heart, segmentation is made even more
difficult by the movement of the organ. As a consequence, the resulting
segmentation is subject to time-dependent shape uncertainty. This uncertainty
in the anatomy propagates through the solution of the inverse problem,
resulting in a reconstruction also affected by uncertainty.

In spite of its acknowledged importance~\cite{clayton2020audit}, uncertainty
quantification in the context of cardiac modelling has emerged only very
recently, see e.g.~\cite{quaglino2018fast,rodriguez2019uncertainty, 
pathmanathan2019comprehensive,corrado2020quantifying}.
In electrocardiography, the most relevant uncertainty to account for
is in the body surface electric recordings. Within a Bayesian framework,
the inverse problem maps a
prior distribution of the pericardial potentials into a posterior distribution
which also accounts for noise in the input data through the
likelihood~\cite{serinagaoglu2005bayesian}. Model uncertainty has been
considered as well. The electric conductivity of the torso is highly
heterogeneous, e.g., lungs, blood masses, interstitial tissue, muscles and
bones significantly differ in terms of conduction, and uncertain, which may
influence the reconstruction as
well~\cite{aboulaich2016stochastic,fikal2019propagation}. In the present
context, shape uncertainty has however received very limited attention thus
far, despite preliminary studies showed a non-negligible impact
on the inverse reconstruction~\cite{tate2018,stoks2019influence}.

From the mathematical perspective, sophisticated tools and theory for the
treatment of shape uncertainties are already available. Besides the fictitious
domain approach considered in \cite{CK07}, one typically distinguishes two
approaches to deal with shape uncertainties: the perturbation method,
see \cite{HSS} and the references therein, which is suitable to treat small
perturbations of the nominal shape and the domain mapping method,
see \cite{xiu,harbrecht2016analysis,GP18,CDE20}.
Here, we focus on the more flexible domain mapping approach. Then, the
computation of quantities of interest, such as expectation and variance of the
potential, gives rise to high dimensional quadrature problems for the random
parameters. Given sufficient parametric regularity, sophisticated sparse
quadrature and quasi-Monte Carlo methods, see e.g.,
\cite{HHPS18,NTW08a,NTW08b,DKL+14,Caf98}, can be applied. If such regularity is
not present, one has to resort to the only slowly converging, Monte Carlo
method.

This work aims at evaluating the impact of space-time uncertainty on the
forward and inverse problems of electrocardiography. To model the shape
uncertainty, we consider a space-time reference geometry and a random
deformation field, that we represent by its Karhunen-Lo\`eve expansion. The
covariance of the random deformation field depends on both space and time and
accounts for the periodicity of the motion of the heart. Furthermore, we
formulate the forward problem as a boundary integral equation by assuming a
constant torso conductivity. We refer to \cite{kress1989linear} for all the
details of this approach and in particular for the applied discretisation.
The boundary integral formulation is particularly suitable for our
purpose, because both input data and observations are confined to a surface,
these are the pericardium and the chest, respectively. In order to accelerate
the computation of quantities of interest, we discretise the spatial problem
with a rapidly converging collocation method, achieving high accuracy already
with a relatively coarse mesh, while the deformation field is approximated by
an efficient low-rank technique.

In order to assess the effect of the shape uncertainty, we estimate the
expectation and the variance of the chest potential (forward problem) and the
pericardial potential (inverse problem) using the anisotropic sparse quadrature
method from \cite{HHPS18}. This approach is validated against the quasi-Monte
Carlo method based on Halton points, see \cite{Caf98}.  The inverse problem,
known to be severely ill-posed, needs an adequate regularisation.  We employ
classic regularisations proposed specifically for our problem of interest,
such as zero order Tikhonov, the first order Tikhonov and total
variation~\cite{karoui2018evaluation}.  The implication of the choice of the
regularisation on the convergence rate of the sparse quadrature is also
analysed. Finally, we explore an $H^{1/2}$ regularisation, see
\cite{OPS10}, a natural option stemming from the boundary integral formulation
of the problem itself. 

The paper is organised as follows. In Section~\ref{sec:setting}, we introduce
the mathematical description of the problem at hand. Section~\ref{sec:BIE}
provides the corresponding boundary integral formulation and the discretisation
of the involved boundary integral operators. The different regularisation methods
for the inverse problem are presented in Section~\ref{sec:inverse}.
Section~\ref{seq:UQ}, dedicated to modelling shape uncertainty, is the
cornerstone of the present work. Finally,
Section~\ref{sec:results} is devoted to the numerical assessment of a
2-dimensional, Cine MRI-derived geometry. Herein, we
validate the sparse quadrature and we quantify the effect of the shape
uncertainty on the forward and inverse problems.

\section{Problem formulation}
\label{sec:setting}

\subsection{The forward problem of electrocardiography}

\begin{figure}[htb]
\centering
\scalebox{0.7}{
\begin{tikzpicture}
\node[inner sep=0pt] (a) {\includegraphics[width=0.4\textwidth]{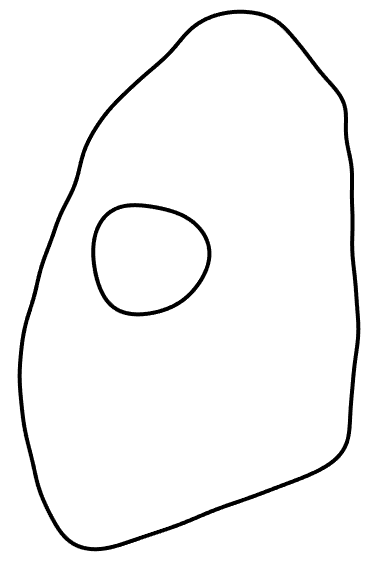}};
\node at ($(a.north)+(1.1,-3.5)$) {\Large$\Sigma(t,\omega)$}
 node at ($(a.north)+(0.5,-6)$) {\Large$D(t,\omega)$}
 node at ($(a.north)+(-1.2,-1.4)$) {\Large$\Gamma$};
\node[anchor=east] (bb) at ($(a.west)-(1cm,0)$) {\includegraphics[width=0.5\textwidth]{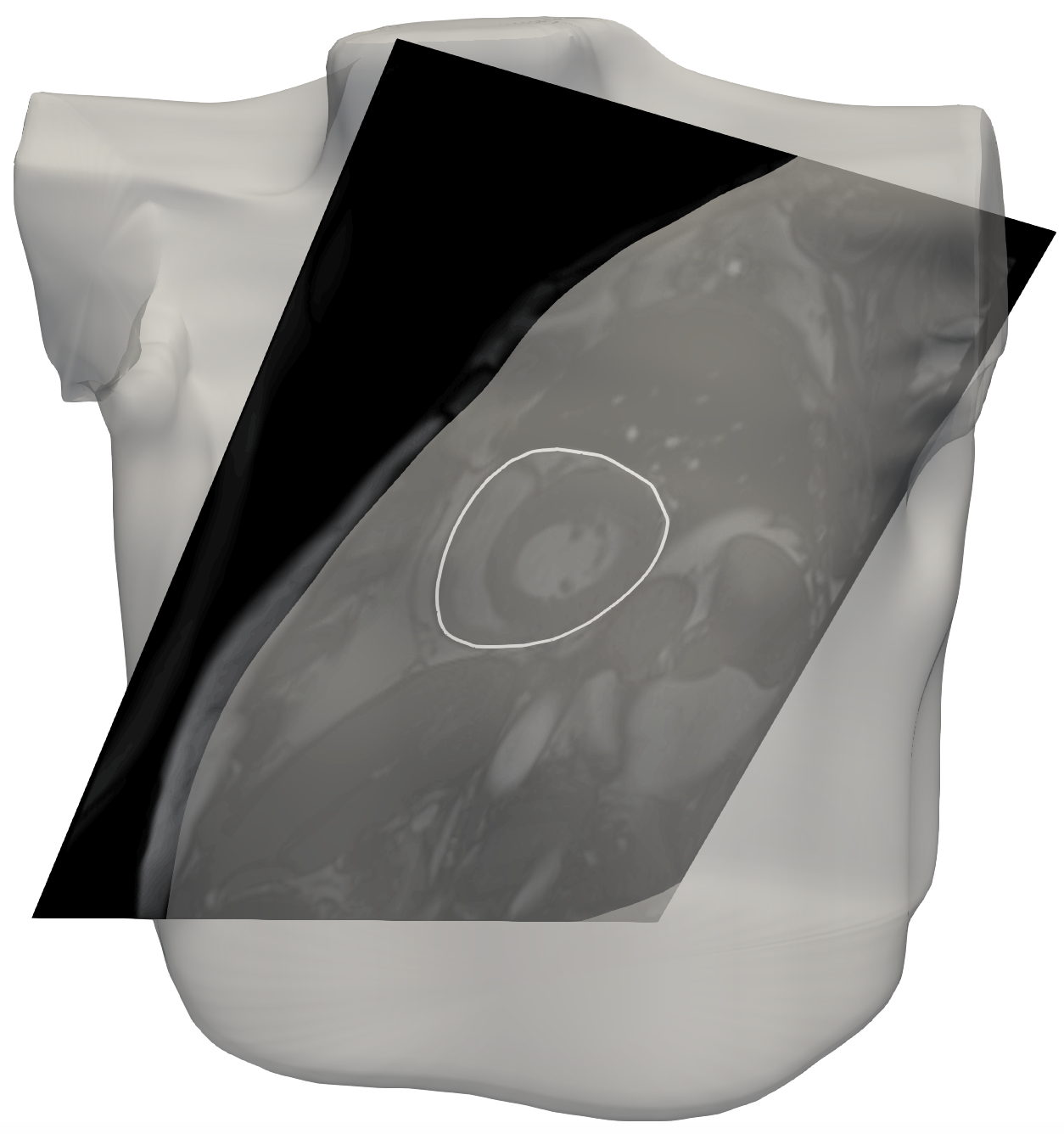}};

\draw[very thick,->] (a.south east) -- node[below] {5\,cm} ($(a.south east)-(1cm,0)$);
\draw[very thick,->] (a.south east) -- node[right] {5\,cm} ($(a.south east)+(0,1cm)$);

\end{tikzpicture}}
\caption{\label{fig:anatomy}2D cross section of the torso reconstructed from
an MRI image.}
\end{figure}
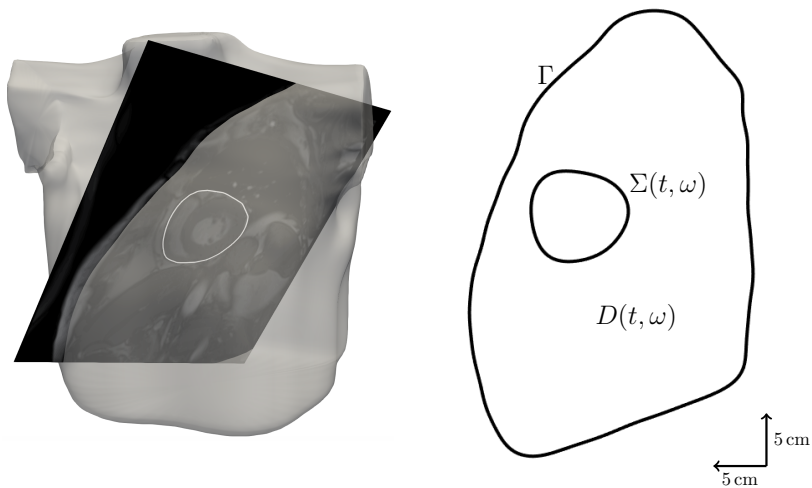

Mathematically, we start from a smooth domain $D\subset\Rbb^d$
representing the torso, whose boundaries are the chest denoted
by $\Gamma$ and the pericardium denoted by $\Sigma$.
In other words, we have $\partial D = \Sigma \cup \Gamma$ and
$\Sigma\cap\Gamma = \emptyset$, see Figure~\ref{fig:anatomy}
for a visualisation of the domain.

The pericardium depends periodically on time and is subject to uncertainty.
To model the time-dependence, let $[0,T]$ be the time interval of interest,
where $T>0$ is the duration of one heartbeat. To model the uncertainty,
let $(\Omega,\mathcal{F},\Pbb)$ denote a complete and separable probability space,
where $\Omega$ is a sample space, $\mathcal{F}\subseteq 2^{\Omega}$ is a
$\sigma$-field and $\Pbb\colon\mathcal{F}\rightarrow[0,1]$
is a probability measure. In what follows, 
we write $\Sigma(t,\omega)$ with \(t\in[0,T)\) and \(\omega\in\Omega\)
to indicate the dependence of the pericardium on time and on the random
parameter. Obviously, the time-dependence and the shape uncertainty of the
heart also affect the torso, which we denote by $D(t,\omega)$,
while we assume the chest $\Gamma$ to be fixed
over time and not being subject to uncertainty.

In the forward problem, given the potential
$u(t)\colon\Sigma(t,\omega)\rightarrow\Rbb$ at the pericardium,
we wish to compute the potential on the chest $\Gamma$. 
Assuming that the torso is a homogeneous volume conductor and that
no current sources are present, the extracellular potential
$y(t,\omega)\colon D(t,\omega)\rightarrow\Rbb$ in the whole torso
satisfies the mixed boundary value problem
\begin{equation}\label{eq:modelProblem}
\left\{
\begin{aligned}
\Delta y(\xbf,t,\omega)&=0, \quad&& \xbf \in D(t,\omega), \\
\frac{\partial y}{\partial{\nbf}_{\xbf}}(\xbf,t,\omega)&= 0,
\quad&& \xbf\in\Gamma,\\
y(\xbf,t,\omega)&= u(\xbf,t), \quad&& \xbf\in\Sigma(t,\omega).
\end{aligned}
\right.
\end{equation}
We remark that in the model the electric conductivity is unit valued,
with no loss of generality since its value, being constant,
does not influence the solution.
Moreover, we assume here that the potential \(u(\xbf,t)\) is given
in spatial coordinates. 
A possibility to guarantee that $u(\xbf,t)$ is well-defined for each realization
of the random parameter is to
assume that it is given with respect to the hold-all domain
\[
\Dcal\isdef\bigcup_{t\in[0,T],\omega\in\Omega}\overline{D(t,\omega)}\times\{t\},
\]
which contains every possible space-time tube.

\subsection{The inverse problem of electrocardiography}

Given the data $y_{\text{d}}(t)\colon\Gamma\rightarrow\Rbb$ on the chest,
the inverse problem corresponding to \eqref{eq:modelProblem} is to
 find the potential $u(t,\omega)\colon\Sigma(t,\omega)\rightarrow\Rbb$
at the pericardium that satisfies
\[
\min_{u(t,\omega)\in L^2(\Sigma(t,\omega))}\|y(t,\omega)-y_{\text{d}}(t)\|_{L^2(\Gamma)}
\quad\mbox{s.t.\ Equation \eqref{eq:modelProblem} holds.}
\]
In the forward problem, if $u(\cdot,t)\in H^{1/2}\big(\Sigma(t,\omega)\big)$, then
the solution satisfies $y(\cdot,t,\omega)\in H^1\big(D(t,\omega)\big)$ and its trace on $\Gamma$ is
in $H^{1/2}(\Gamma)$. 
This gives rise to the solution operator
\[
\Acal(t,\omega)\colon H^{1/2}\big(\Sigma(t,\omega)\big)
\to H^{1/2}(\Gamma).
\]
It is possible to show, see \cite{colli1985inverse},
that \(\Acal(t,\omega)\) is injective and continuous for \(t\in[0,T]\) and
almost every \(\omega\in\Omega\). In the inverse problem,
therefore, if $y_{\text{d}}(\cdot,t)\in H^{1/2}(\Gamma)$ then there exists
a unique minimum $u(\cdot,t,\omega)\in H^{1/2}\big(\Sigma(t,\omega)\big)$
such that $y(t,\omega)=y_{\text{d}}(t)$ at $\Gamma$.
The minimum $u(t,\omega)$, however, is not stable with respect to perturbations.
A common practice is to introduce a regularisation into the problem to
recover stability, as follows:
\[
\min_{u(t,\omega)\in L^2(\Sigma(t,\omega))}\left\{ \frac{1}{2} \|
\Acal(t,\omega)u(t,\omega) - y_{\text{d}}(t) \|^2_{L^2(\Gamma)}
+ \frac{\lambda}{2}\mathcal{R}\big(u(t,\omega)\big)\right\}.
\]
Herein, $\lambda>0$ is the regularisation parameter and $\mathcal{R}$ is the
regularisation functional, for which several choices are possible. 
The optimum can be explicitly computed according to
\[
\langle \Acal(t,\omega)u(t,\omega) - y_{\text{d}}(t),\Acal(t,\omega) \delta u\rangle_{L^2(\Gamma)}
+ \lambda \langle \mathcal{R}'\big(u(t,\omega)\big), \delta u\rangle = 0, \quad \forall\delta u \in L^2(\Sigma(t,\omega)).
\]
Later on
we shall, discuss different options for the regularisation.

\section{Boundary integral formulation}
\label{sec:BIE}
\subsection{Boundary integral operators}\label{sec:BIOs}
In the following exposition, for the sake of an easier notation,
we consider $t$ and $\omega$ to be fixed.
For a comprehensive exposition of the boundary integral approach, we refer the
reader to \cite{kress1989linear}. For a
a better distinction of quantities that are defined with respect to
the boundaries and those which are given on the domain, we introduce 
for \(\Phi\in\{\Sigma,\Gamma,\partial D\}\) the
trace operators
\begin{equation*}
\gamma_{0,\Phi}^{\operatorname{int}}\colon H^1(D)\to H^{1/2}(\Phi)
\quad\mbox{and}\quad
\gamma_{1,\Phi}^{\operatorname{int}}\colon H^1(D)\to H^{-1/2}(\Phi).
\end{equation*}
Then, we may write the forward problem according to
\begin{equation*}
\begin{cases}
\Delta y(\xbf)=0,& \xbf\in D, \\
\gamma_{1,\Gamma}^{\operatorname{int}}y(\xbf)= 0,& \xbf\in\Gamma, \\
\gamma_{0,\Sigma}^{\operatorname{int}}y(\xbf)= u(\xbf), & \xbf\in\Sigma.
\end{cases}
\end{equation*}
For $\xbf\in D$, the potential \(y(\xbf)\) is given by the representation formula
\begin{equation}
\label{eq:repForm}
y(\xbf) = \int_{\partial D}G(\xbf,\xbf^\prime)(\gamma_{1,\partial D}^{\operatorname{int}}y)(\xbf^\prime)\dd\sigma_{\xbf^\prime}
-\int_{\partial D}\frac{\partial G}{\partial\nbf_{\xbf^\prime}}(\xbf,\xbf^\prime)(\gamma_{0,\partial D}^{\operatorname{int}}y)(\xbf^\prime)\dd\sigma_{\xbf^\prime},
\end{equation}
where
\begin{equation}\label{fundSol}
G(\xbf,\xbf')\isdef
\begin{cases}
-\dfrac{1}{2\pi}\log\|\xbf - \xbf'\|_2,& d=2,\\[1ex]
\phantom{-}\dfrac{1}{4\pi}\dfrac{1}{\|\xbf-\xbf'\|_2},& d=3,
\end{cases}
\end{equation}
denotes the fundamental solution for the Laplacian. 
Introducing the single layer operator 
\[
\Vcal\colon H^{-1/2}(\partial D)\to H^{1/2}(\partial D),\quad
(\Vcal \rho)(\xbf)\isdef\int_{\partial D}G(\xbf,\xbf^\prime)\rho(\xbf^\prime)\dd\sigma_{\xbf^\prime},
\] 
the double layer operator
\[
\Kcal\colon H^{1/2}(\partial D)\to H^{1/2}(\partial D),\quad
(\Kcal \rho)(\xbf)\isdef\int_{\partial D}\frac{\partial G}
{\partial\nbf_{\xbf^\prime}}(\xbf,\xbf^\prime)\rho(\xbf^\prime)\dd\sigma_{\xbf^\prime}
\]
and taking the Dirichlet trace, Equation~\eqref{eq:repForm}
yields the Dirichlet-to-Neumann map
\begin{equation}\label{eq:BIE}
\Vcal(\gamma_{1,\partial D}^{\operatorname{int}}y)=
\bigg(\frac 1 2 I + \Kcal\bigg)(\gamma_{0,\partial D}^{\operatorname{int}}y).
\end{equation}
Next, introducing for \(\Phi,\Psi\in\{\Sigma,\Gamma\}\), the restricted operators
\[
\Vcal_{\Phi\Psi}\colon H^{-1/2}(\Psi)\to H^{1/2}(\Phi),
\quad \Phi\ni \xbf \mapsto
(\Vcal_{\Phi\Psi}\rho)(\xbf)\isdef \int_{\Psi}G(\xbf,\xbf^{\prime})\rho(\xbf^{\prime})\dd\sigma_{\xbf^{\prime}}
\]
and
\[
\Kcal_{\Phi\Psi}\colon H^{1/2}(\Psi)\to H^{1/2}(\Phi),
\quad \Phi\ni \xbf \mapsto
(\Kcal_{\Phi\Psi}\rho)(\xbf)\isdef \int_{\Psi}\frac{\partial G}
{\partial\nbf_{\xbf^\prime}}(\xbf,\xbf^{\prime})\rho(\xbf^{\prime})
\dd\sigma_{\xbf^{\prime}},
\]
we can split up Equation~\eqref{eq:BIE} into the system
\begin{equation*}
\begin{bmatrix}
\mathcal{V}_{\Sigma\Sigma} & \mathcal{V}_{\Sigma\Gamma}\\
\mathcal{V}_{\Gamma\Sigma} & \mathcal{V}_{\Gamma\Gamma}
\end{bmatrix}
\begin{bmatrix}
\gamma_{1,\Sigma}^{\text{int}}y \\
\gamma_{1,\Gamma}^{\text{int}}y
\end{bmatrix}=
\begin{bmatrix}
\frac{1}{2}I+\mathcal{K}_{\Sigma\Sigma} & \mathcal{K}_{\Sigma\Gamma}\\
\mathcal{K}_{\Gamma\Sigma} & \frac{1}{2}I+\mathcal{K}_{\Gamma\Gamma} 
\end{bmatrix}
\begin{bmatrix}
\gamma_{0,\Sigma}^{\text{int}}y \\
\gamma_{0,\Gamma}^{\text{int}}y
\end{bmatrix}.
\end{equation*}
Rearranging this system in order to move the unknowns
$\gamma_{1,\Sigma}^\text{int}y$ and $\gamma_{0,\Gamma}^\text{int}y$
to the left side and the data $\gamma_{0,\Sigma}^\text{int}y=u$
and $\gamma_{1,\Gamma}^\text{int}y=0$ to the right side,
we arrive at the system boundary integral equations
\begin{equation}\label{BoundaryIntEq}
\begin{bmatrix}
\mathcal{V}_{\Sigma\Sigma}& -\mathcal{K}_{\Sigma\Gamma}\\
-\mathcal{V}_{\Gamma\Sigma} & \frac{1}{2}I+\mathcal{K}_{\Gamma\Gamma}
\end{bmatrix}
\begin{bmatrix}
\gamma_{1,\Sigma}^{\text{int}}y \\
\gamma_{0,\Gamma}^{\text{int}}y
\end{bmatrix}=
\begin{bmatrix}
\frac{1}{2}I+\mathcal{K}_{\Sigma\Sigma} & -\mathcal{V}_{\Sigma\Gamma}\\
-\mathcal{K}_{\Gamma\Sigma} & \mathcal{V}_{\Gamma\Gamma} 
\end{bmatrix}
\begin{bmatrix}
u \\ 0
\end{bmatrix}.
\end{equation}
As \(\mathcal{V}_{\Sigma\Sigma}\) is an elliptic operator,
considering the Schur complement 
\[
\mathcal{S}\colon H^{1/2}(\Gamma)\to H^{1/2}(\Gamma),
\quad \mathcal{S}\rho \isdef
\bigg(\frac{1}{2}I+\mathcal{K}_{\Gamma\Gamma}-\mathcal{V}_{\Gamma\Sigma}\mathcal{V}_{\Sigma\Sigma}^{-1}\mathcal{K}_{\Sigma\Gamma}
\bigg)\rho
\]
yields
the operator equation
\[
\mathcal{S}\gamma_{0,\Gamma}^{\text{int}}y
=\bigg(\frac 1 2\mathcal{V}_{\Gamma\Sigma}\mathcal{V}_{\Sigma\Sigma}^{-1}+
\mathcal{V}_{\Gamma\Sigma}\mathcal{V}_{\Sigma\Sigma}^{-1}
\mathcal{K}_{\Sigma\Sigma}-\mathcal{K}_{\Gamma\Sigma}\bigg)u
\]
for the desired potential \(\gamma_{0,\Gamma}^{\text{int}}y\) on the chest.
In particular, the solution operator \(\mathcal{A}\colon H^{1/2}(\Sigma)
\to H^{1/2}(\Gamma)\) is given by
\begin{equation}\label{eq:SolOp}
\mathcal{A}u=
\mathcal{S}^{-1}\bigg(\frac 1 2\mathcal{V}_{\Gamma\Sigma}
\mathcal{V}_{\Sigma\Sigma}^{-1}+
\mathcal{V}_{\Gamma\Sigma}\mathcal{V}_{\Sigma\Sigma}^{-1}
\mathcal{K}_{\Sigma\Sigma}-\mathcal{K}_{\Gamma\Sigma}\bigg)u
=\gamma_{0,\Gamma}^{\text{int}}y,
\end{equation}
while the Poincar\'e-Steklov operator
\(\mathcal{B}\colon H^{1/2}(\Sigma)\to H^{-1/2}(\Sigma)\) is given by
\begin{equation}\label{eq:PSop}
\mathcal{B}u=\mathcal{V}_{\Sigma\Sigma}^{-1}
\bigg(\frac{1}{2}I+\mathcal{K}_{\Sigma\Sigma}
+\mathcal{K}_{\Sigma\Gamma}\mathcal{A}\bigg)u
=\gamma_{1,\Sigma}^{\text{int}}y.
\end{equation}

\subsection{Numerical discretisation}
Our goal is to discretise the boundary integral equations derived previously.
In order to not having to deal with non-constant coefficients,
we shall solve the equations in the
spatial coordinate frame, i.e., for each tuple \((t,\omega)\),
assemble the boundary integral
operators on \(\Sigma(t,\omega)\) and \(\Gamma\).
To simplify the notation, in the following quantities we omit the dependency
of $\Sigma$ on $t$ and $\omega$.
As the focus of this work is on uncertainty quantification, rather than on
the solution of the spatial problem, we restrict ourselves to
the simplified 2D anatomy of Figure~\ref{fig:anatomy}. Hence, from \eqref{fundSol},
we have
\[
G(\xbf,\xbf^\prime)=-\frac{1}{2\pi}\log\|\xbf-\xbf^\prime\|_2
\]
and
\[
\frac{\partial G}{\partial\nbf_{\xbf^\prime}}(\xbf,\xbf^\prime)=-\frac{1}{2\pi}\nabla_{\xbf^\prime}\log\|\xbf-\xbf^\prime\|_2\cdot\nbf_{\xbf^\prime}=\frac{1}{2\pi}\frac{\langle \xbf-\xbf^\prime,\nbf_{\xbf^\prime}\rangle_2}{\|\xbf-\xbf^\prime\|_2^2}.
\]
Given that the domain \(D(t,\omega)\)
exhibits \(C^2\)-boundaries, which are described by the parameterisations
$\gamma_{\Sigma}\colon[0,1)\rightarrow\Sigma(t,\omega)$ and $\gamma_{\Gamma}\colon [0,1)\rightarrow\Gamma$, we can employ the collocation method from \cite{kress1989linear}. This method is based on the trapezoidal rule and an
appropriate desingularisation of the kernel functions. In view of the
Euler-Maclaurin formula, it converges
quadratically for \(C^2\)-boundaries and exponentially for smooth boundaries.
More precisely,
given that the solution to \eqref{BoundaryIntEq} satisfies 
\(y\in C^k(\partial D)\), there holds
\[
\|y-y_n\|_{L^\infty(\partial D)}\leq C n^{-k}\|y\|_{C^k(\partial D)}
\]
for some \(C>0\),
where \(y_n\) is obtained from the collocation method with \(n=2j\) points for some 
\(j\in\mathbb{N}\). We refer to \cite{kress1989linear}
for all the details.
For the sake of completeness, we briefly recall the collocation method in Appendix~\ref{sec:Nystrom}.

For $\Phi\in\{\Sigma,\Gamma\}$, we consider the $n_{\Phi}$-points discretisations
$s_{i,\Phi}=i/n_{\Phi}$ for $i=0,\ldots,n_{\Phi}-1$.
Moreover, we denote by
$\Ibf_{n_{\Phi}}\in\mathbb{R}^{n_{\Phi}\times n_{\Phi}}$ the $n_{\Phi}\times n_{\Phi}$
identity matrix. The discretisation results in the block linear system
\begin{equation}\label{LinSys}
\begin{bmatrix}
\Vbf_{\Sigma\Sigma}& -\Kbf_{\Sigma\Gamma}\\
-\Vbf_{\Gamma\Sigma} & \frac{1}{2}\Ibf_{n_{\Gamma}}+\Kbf_{\Gamma\Gamma}
\end{bmatrix}
\begin{bmatrix}
\tilde{\boldsymbol{\rho}}_{1,\Sigma} \\
\boldsymbol{\rho}_{0,\Gamma}
\end{bmatrix}=
\begin{bmatrix}
\frac{1}{2}\Ibf_{n_{\Sigma}}+\Kbf_{\Sigma\Sigma} & -\Vbf_{\Sigma\Gamma}\\
-\Kbf_{\Gamma\Sigma} &  \Vbf_{\Gamma\Gamma} 
\end{bmatrix}
\begin{bmatrix}
{\mathbf {u}} \\ {\mathbf{0}_{n_{\Gamma}}}
\end{bmatrix},
\end{equation}
cp.\ \eqref{BoundaryIntEq}. Herein
\[
\boldsymbol{\rho}_{0,\Gamma}\isdef\big[y\big(\gamma_{\Gamma}(s_{i,\Gamma})\big)\big]_{i=0}^{n_{\Gamma}-1}\in\mathbb{R}^{n_{\Gamma}}
\quad\text{and}\quad
\tilde{\boldsymbol{\rho}}_{1,\Sigma}\isdef\left[\frac{\partial y}{\partial\nbf_{s_i,\Sigma}}\big(\gamma_{\Sigma}(s_{i,\Sigma})\big)\big\|\gamma_{\Sigma}^{\prime}(s_{i,\Sigma})\big\|_2\right]_{i=0}^{n_{\Sigma}-1}\in\mathbb{R}^{n_{\Sigma}}.
\]
Furthermore, we set
\[
\mathbf{u}\isdef\big[u\big(\gamma_{\Sigma}(s_{i,\Sigma})\big)\big]_{i=0}^{n_{\Sigma}-1}\in\mathbb{R}^{n_{\Sigma}}
\]
and $\mathbf{0}_{n_{\Gamma}}\in\mathbb{R}^{n_{\Gamma}}$ is
the $n_{\Gamma}$-dimensional vector of zeros. 

\section{Solution of the inverse problem}
\label{sec:inverse}

In the inverse problem, the goal is to reconstruct the potential
$u$ at the pericardium $\Sigma(t,\omega)$ given noisy data measurements
$y_{\text{d}}\in H^{1/2}(\Gamma)$ on the chest.
In the following exposition, for the sake of an easier notation,
we indicate the dependency of $\Sigma$ on the tuple $(t,\omega)$ only when
specifying function spaces.
In particular, we consider the solution operator
$\mathcal{A}\colon H^{1/2}\big(\Sigma(t,\omega)\big)\rightarrow H^{1/2}(\Gamma)$,
cp.\ \eqref{eq:SolOp}, and the Poincar\'e-Steklov operator $\mathcal{B}\colon H^{1/2}\big(\Sigma(t,\omega)\big)\rightarrow H^{-1/2}\big(\Sigma(t,\omega)\big)$,
cp.\ \eqref{eq:PSop}.
In the discrete setting, we define the vectors
\[
\mathbf{u}\isdef\big[u\big(\gamma_{\Sigma}(s_{i,\Sigma})\big)\big]_{i=0}^{n_{\Sigma}-1}\in\mathbb{R}^{n_{\Sigma}},\quad \mathbf{y}_{\text{d}}\isdef
\big[y_{\text{d}}\big(\gamma_{\Gamma}(s_{i,\Gamma})\big)\big]_{i=0}^{n_{\Gamma}-1}\in\mathbb{R}^{n_{\Gamma}}
\]
and
\[
\boldsymbol{\rho}_{0,\Gamma}\isdef\big[y\big(\gamma_{\Gamma}(s_{i,\Gamma})\big)\big]_{i=0}^{n_{\Gamma}-1}\in\mathbb{R}^{n_{\Gamma}},\quad \boldsymbol{\rho}_{1,\Sigma}\isdef\left[\frac{\partial y}{\partial\mathbf{n}_{s_i,\Sigma}}\big(\gamma_{\Sigma}(s_{i,\Sigma})\big)\right]_{i=0}^{n_{\Sigma}-1}\in\mathbb{R}^{n_{\Sigma}}.
\]
From \eqref{LinSys} we deduce that there exist a matrix $\Abf\in\mathbb{R}^{n_{\Gamma}\times n_{\Sigma}}$ and a matrix $\Bbf\in\mathbb{R}^{n_{\Sigma}\times n_{\Sigma}}$ such that
\[
\boldsymbol{\rho}_{0,\Gamma}=\Abf\mathbf{u} \quad\text{and}\quad \boldsymbol{\rho}_{1,\Sigma}=\Bbf\mathbf{u}.
\]
In this section we consider inverse problem solutions of the form
\begin{equation}\label{contInv}
u=\argmin_{v\in L^2(\Sigma(t,\omega))}\left\{\frac{1}{2}\|\mathcal{A}v-y_{\text{d}}\|^2_{L^2(\Gamma)}+\frac{\lambda}{2}\mathcal{R}(v)\right\},
\end{equation}
where $\lambda>0$ is the regularisation parameter and $\mathcal{R}$ is the regularisation functional.
We discretise \eqref{contInv} employing, for $\Phi\in\{\Sigma,\Gamma\}$, the trapezoidal rule
\[
\langle v,w\rangle_{L^2(\Phi)}\approx{\bf v}^\intercal\Sbf_{\Phi}{\bf w}\quad\text{with}\quad{\bf v}\isdef[v\big(\gamma_{\Phi}(s_{i,\Phi})\big)]_{i=0}^{n_{\Phi}-1},\ {\bf w}\isdef[w\big(\gamma_{\Phi}(s_{i,\Phi})\big)]_{i=0}^{n_{\Phi}-1},
\]
where $\Sbf_{\Phi}\in\mathbb{R}^{n_{\Phi}\times n_{\Phi}}$ is the
diagonal mass matrix with entries
\[
(\Sbf_{\Phi})_{i,i}\isdef\frac{\|\gamma_{\Phi}^\prime(s_{i,\Phi})\|_2}{n_{\Phi}}
\]
for $i=0,\ldots,n_{\Phi}-1$. The accuracy of this approximation is consistent
with the one obtained by the collocation method, as can easily be inferred
from the Euler-Maclaurin formula. 

We arrive at the discrete formulation
\begin{equation}\label{discrInv}
\mathbf{u}=\argmin_{\mathbf{v}\in\mathbb{R}^{n_{\Sigma}}}
\left\{\frac{1}{2}(\Abf\mathbf{v}-\mathbf{y}_{\text{d}})^\intercal\Sbf_{\Gamma}
(\Abf\mathbf{v}-\mathbf{y}_{\text{d}})+\frac{\lambda}{2}R(\mathbf{v})\right\},
\end{equation}
where $R$ is the discrete regularisation term. The optimisation of the objective
function leads to the equation
\begin{equation}\label{discrOpt}
\Abf^\intercal\Sbf_{\Gamma}(\Abf\mathbf{u}-\mathbf{y}_{\text{d}})+\frac{\lambda}{2}R^\prime(\mathbf{u})=0
\end{equation}
for the minimum $\mathbf{u}$. In the following subsections we present four
strategies to solve the inverse problem based on different regularisations.

\subsection{Zero order Tikhonov regularisation}

The zero order Tikhonov regularisation penalises the $L^2$-norm of the
Dirichlet data at $\Sigma(t,\omega)$, see \cite{mueller2012linear},
and the regularisation functional in \eqref{contInv} reads
\[
\mathcal{R}(v)=\|v\|^2_{L^2(\Sigma(t,\omega))}.
\]
The discrete regularisation term in \eqref{discrInv} is thus
\[
R(\mathbf{v})=\mathbf{v}^\intercal\mathbf{S}_{\Sigma}\mathbf{v}
\]
and in \eqref{discrOpt} we have 
\[
R^\prime(\mathbf{u})=2\mathbf{S}_{\Sigma}\mathbf{u}.
\]

\subsection{First order Tikhonov regularisation}

The first order Tikhonov regularisation penalises the $L^2$-norm of the Neumann
data at $\Sigma(t,\omega)$, see \cite{mueller2012linear}, and the regularisation
functional in \eqref{contInv} is
\[
\mathcal{R}(v)=\|\mathcal{B}v\|^2_{L^2(\Sigma(t,\omega))}.
\]
Then the discrete regularisation term in \eqref{discrInv} is
\[
R(\mathbf{v})=\mathbf{v}^\intercal\mathbf{B}^\intercal
\mathbf{S}_{\Sigma}\mathbf{B}\mathbf{v}
\]
and in \eqref{discrOpt} we have 
\[
R^\prime(\mathbf{u})=2\mathbf{B}^\intercal
\mathbf{S}_{\Sigma}\mathbf{B}\mathbf{u}.\]

\subsection{\texorpdfstring{$H^{1/2}$}{H1/2} regularisation}

When penalising the $H^{1/2}$-norm of the Dirichlet data at $\Sigma(t,\omega)$,
the regularisation functional in \eqref{contInv} is given by
\[
\mathcal{R}(v)=\|v\|^2_{H^{1/2}(\Sigma(t,\omega))}
=\langle\mathcal{B}v,v\rangle_{L^2(\Sigma(t,\omega))}.
\]
The resulting discrete regularisation term in \eqref{discrInv} is 
\[
R(\mathbf{v})=\mathbf{v}^\intercal\mathbf{B}^\intercal
\mathbf{S}_{\Sigma}\mathbf{v}
\]
and its Fr\'echet derivative in \eqref{discrOpt} is
\[
R^\prime(\mathbf{u})=\mathbf{S}_{\Sigma}(\mathbf{B}^\intercal
+\mathbf{B})\mathbf{u}.
\]

\subsection{Total Variation regularisation}

In the Total Variation regularisation we penalise the $L^1$-norm of
the Neumann data at $\Sigma(t,\omega)$, i.e., the regularisation functional
in \eqref{contInv} is
\[
\mathcal{R}(v)=\|\mathcal{B}v\|_{L^1(\Sigma(t,\omega))}.
\]
As the $L^1$-norm is non-differentiable, it is common to employ
the approximation
\[
\|\mathcal{B}v\|_{L^1(\Sigma(t,\omega))}=\int_{\Sigma(t,\omega)}|
\mathcal{B}v(\xbf)|\text{d}\sigma_{{\bf x}}\approx\int_{\Sigma(t,\omega)}
\sqrt{(\mathcal{B}v(\boldsymbol{x}))^2+\beta}\ \text{d}\sigma_{{\bf x}},
\]
where $\beta>0$ is a small constant (typically $\beta=10^{-5}$).
This introduces a non-linearity in the optimality condition.
This difficulty is handled in \cite{karoui2018evaluation}
by using the zero order Tikhonov solution $u_0$ in the non-linear term.
This approach is equivalent to set
\[
\mathcal{R}(v)=\int_{\Sigma(t,\omega)}
\Big(2\sqrt{\big(\mathcal{B}u_0(\mathbf{x})\big)^2+\beta}\Big)^{-1}\big(\mathcal{B}v(\mathbf{x})\big)^2\dd\sigma_{\mathbf{x}}
\]
in \eqref{contInv}.
The discretisation requires the definition of the diagonal matrix
$\mathbf{W}_{\beta}(\mathbf{u}_0)\in\mathbb{R}^{n_{\Sigma}\times n_{\Sigma}}$
with entries
\[
\big(\Wbf_{\beta}(\mathbf{u}_0)\big)_{i,i}
\isdef\Big(2\sqrt{(\Bbf\mathbf{u}_0)_i^2+\beta}\Big)^{-1}
\]
for $i=0,\ldots,n_{\Sigma}-1$, where $\mathbf{u}_0$ is the discrete
zero order Tikhonov solution.
The discrete regularisation term in \eqref{discrInv} then reads
\[
R(\mathbf{v})=\mathbf{v}^\intercal\mathbf{B}^\intercal\mathbf{W}_{\beta}(\mathbf{u}_0)\mathbf{S}_{\Sigma}\mathbf{B}\mathbf{v},
\]
while in \eqref{discrOpt} it results in
\[
R^\prime(\mathbf{u})=2\mathbf{B}^\intercal
\mathbf{W}_{\beta}(\mathbf{u}_0)\mathbf{S}_{\Sigma}\mathbf{B}\mathbf{u}.
\]

\section{Uncertainty Quantification}
\label{seq:UQ}

\subsection{Modeling time-dependent shape uncertainty}
\label{sec:uncertainty}

In order to assess the shape uncertainty of the pericardium $\Sigma(t,\omega)$,
we assume the existence of a reference domain \(D_\refd\subset\mathbb{R}^d\),
with boundaries $\Sigmaref$ and \(\Gamma\),
and a random deformation field 
\[\chibf\colon D_\refd
\times[0,T]\times\Omega\rightarrow\mathbb{R}^d\] such that
\[
D(t,\omega)=\chibf(D_\refd,t,\omega)=\chibf\big(D_\refd(t),\omega\big),
\]
where 
\[
D_\refd(t)=\Ebb[\chibf](D_\refd,t).
\]
Herein, the expectation is given in terms of the Bochner integral
\[
\Ebb[\chibf](\xmat,t)\isdef\int_\Omega \chibf(\xmat,t,\omega)\dd\Pbb.
\]
Introducing the space-time cylinder
\(
\Qref\isdef D_\refd\times(0,T)
\) with coordinates 
\[
(\xmat,t)\defis \zmat\in \Qref,
\] we can interpret \(\chibf\)
as a random deformation field \(\chibf(\zmat,\omega)\) acting on the 
space-time cylinder. A similar model has originally been introduced in \cite{harbrecht2016analysis} for the stationary case.
Assuming \[\chibf\in L^2\big(\Omega; C^1(Q_\refd;\mathbb{R}^d)\big),\]
it is possible to expand \(\chibf(\zmat,\omega)\) in a Karhunen-Lo\`eve
expansion, cf.\ \cite{Loe77}, according to
\[
\chibf(\zmat,\omega)=\Ebb[\chibf](\zmat)+\sum_{k=1}^\infty
\sqrt{\lambda_k}\chibf_k(\zmat)Y_k(\omega).
\]
Herein, \((\lambda_k,\chibf_k)\) are the eigen-pairs of the integral
operator defined by the
matrix valued covariance function
\begin{align*}
\Cov[\chibf]\colon \Qref\times \Qref\to\mathbb{R}^{d\times d},\
\Cov[\chibf](\zmat,\zmat')\isdef\int_{\Omega}
\big(\chibf(\zmat,\omega)-\Ebb[\chibf](\zmat)\big)
\big(\chibf(\zmat',\omega)-\Ebb[\chibf](\zmat')\big)^\intercal\dd\Pbb,\nonumber
\end{align*}
while \(\{Y_k\}\) is a family of uncorrelated random variables
with normalised variance.

To avoid distortion of the space-time tubes \(\chibf(\Qref,\omega)\),
\(\omega\in\Omega\), and hence to guarantee well-posedness of the boundary
value problem at hand, cf.\ \eqref{eq:modelProblem},
we impose the uniformity condition
\begin{equation*}
\|\chibf(\cdot,\omega)\|_{C^1(\overline{\Qref};\mathbb{R}^d)},
\|\chibf^{-1}(\cdot,\omega)\|_{C^1(\overline{Q(\omega)};\mathbb{R}^d)}
< C_{\mathrm{uni}}
\end{equation*}
for some constant \(C_{\mathrm{uni}}>0\) uniformly in
\(\omega\in\Omega\).\footnote{%
For the collocation method, we require also \(\chibf(t,\omega)
\in C^2(\overline{D_\refd})\) for
\(t\in[0,T],\omega\in\Omega\).}
As a consequence, the random variables' \(\{Y_k\}\) ranges are bounded.
Therefore, without loss of
generality we may assume \(Y_k\colon\Omega\to[-1,1]\), \(k=1,2,\ldots\).
In particular, we obtain
\[
\Sigma(t,\omega)=\chibf(\Sigmaref,t,\omega)
\quad\text{and}\quad\Gamma=\chibf(\Gamma,t,\omega).
\]
Replacing the random variables by their image \([-1,1]^\Nbb\),
we arrive at the parametrised Karhunen-Lo\`eve expansion
\begin{equation*}
\chibf(\zmat,\rparam)=\Ebb[\chibf](\zmat)+\sum_{k=1}^\infty
\sqrt{\lambda_k}\chibf_k(\zmat)\rparamscal_k,\quad
\rparam\in[-1,1]^\Nbb.
\end{equation*}
For the numerical computation of the Karhunen-Lo\`eve expansion it is sufficient
to know the expectation
and the covariance at \(\Sigma_\refd\times[0,T]\), see \cite{DHJM20}. In our concrete case,
it is even sufficient to only compute the
Karhunen-Lo\`eve expansion in the space-time collocation points. To this end, we employ
the pivoted Cholesky decomposition, see \cite{HPS12,HPS15,harbrecht2016analysis}.
This results in a finite rank Karhunen-Lo\`eve expansion
\begin{equation}\label{eq:finiteKL}
\chibf(\zmat_i,\rparam)=\Ebb[\chibf](\zmat_i)+\sum_{k=1}^K
\sqrt{\lambda_k}\chibf_k(\zmat_i)\rparamscal_k,\quad
\rparam\in[-1,1]^K,
\end{equation}
where $\zmat_i$ are the space-time collocation points on $\Sigma_{\text{ref}}\times[0,T]$ and $K$ is the dimension of the random parameter space.

\subsection{Computation of quantities of interest}
We make the common assumption that the random variables are even independent and 
uniformly distributed. Then, the push-forward measure is given by the product
density
\[
\rho({\rparam})\isdef\prod_{k=1}^K\rho_k(\rparamscal_k)\quad\text{with}\quad
\rho_k\equiv\frac 1 2.
\]
Hence, we can now express quantities of interest, such as the moments of the
potential on the chest, by means of an integral with respect to \([-1,1]^K\).
It holds for the moments that
\[
\mathcal{M}_m[y]({\bf x},t)\isdef
\int_{\Omega}y({\bf x},t,\omega)^m\dd\mathbb{P}
=\int_{[-1,1]^K}y(\xbf,t,\rparam)^m\rho(\rparam)\dd\rparam,
\quad m\in\Nbb^*.
\]
Especially, the expectation and the variance are given by
\[
\mathbb{E}[y]({\bf x},t)=\mathcal{M}_1[y]({\bf x},t)
\quad\text{and}\quad \mathbb{V}[y]({\bf x},t)
=\mathcal{M}_2[y]({\bf x},t)-\big(\mathcal{M}_1[y]({\bf x},t)\big)^2.
\]
In practice, the moments need to be approximated by a quadrature rule,
i.e.,
\[
\mathcal{M}_m[y]({\bf x},t)\approx\sum_{i=1}^N w_i y({\bf x},t,\rparam_i)^m,
\]
where $\rparam_i\in[-1,1]^K$ for $i=1,\ldots,N$ are the quadrature points
and $w_i\in\mathbb{R}$ are the corresponding weights.
Due to the space-time modelling, the random parameter space is very
high-dimensional, and appropriate quadrature rules need to be  employed.
In this article, we employ the anisotropic sparse quadrature from \cite{HHPS18},
which is based on the Gauss-Legendre quadrature rule. For the sparse quadrature
to converge, the quantity of interest needs to be regular with respect to
the random parameter $\rparam\in[-1,1]^K$. For the forward problem,
such results are available, we refer to \cite{harbrecht2016analysis}. For
the Dirichlet control problem with an affine random diffusion coefficient,
such results have also been derived in \cite{KS13}.

\section{Numerical experiments}
\label{sec:results}

In this section we present numerical experiments to assess
the validity of the approach. To obtain a space-time reference torso anatomy,
we have segmented cardiac magnetic resonance (CMR) imaging data previously acquired~\cite{Potse2014}. The temporal image stack was obtained from  a cine
ECG-triggered segmented steady-state free precession sequence in
mid-ventricular short-axis orientation. Slice thickness and voxel resolution
were \SI{8}{\mm} and \SI{1.36}{\mm}, respectively.
The 25-phase temporal stack covered the whole heart beat, with an
RR interval on the surface ECG of $T=\SI{690}{\ms}$).
Images were reordered such that the first image
at $t=\SI{0}{\ms}$ corresponded to the diastole, defined by the maximal left
ventricular cavity volume. The systole, defined by the minimum cavity volume,
occurred at $t=\SI{270}{\ms}$. The segmentation was performed by manual contour
tracing of the epicardium for each image. In order to end up with a smooth
computational domain, we performed a least-squares fit of the contours using
a truncated Fourier series with a threshold of $10^{-3}$ in the relative
root-mean-square error.  Finally, we interpolated the extracted shapes to get a
pericardial representation at 50 time instants. The chest was also previously
segmented from ultra-fast gradient-echo ``VIBE'' images in axial, coronal, and
sagittal orientations to produce a smooth 3-dimensional closed surface modelled
in Blender\footnote{\texttt{https://www.blender.org}}. 
As shown in Figure~\ref{fig:anatomy}, the contour of the chest was eventually
obtained by intersecting the chest surface
with the orientation plane of the pericardium.

\begin{figure}[t]
\includegraphics[width=\textwidth, keepaspectratio]{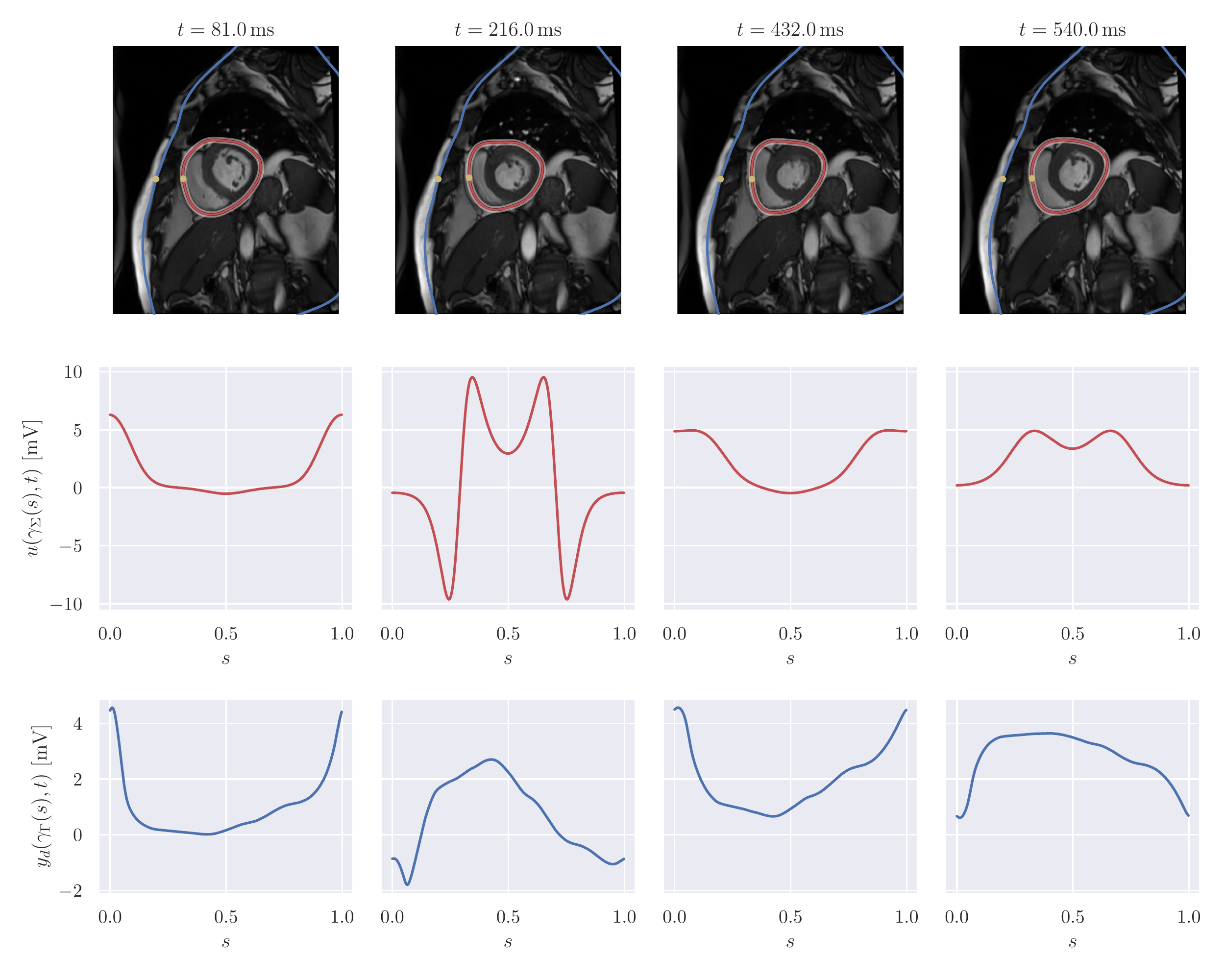}
\caption{Geometry and input data.
First row: CMR images with superimposed segmented chest (blue) and
time-dependent pericardium (red).  The shaded region around the pericardium
represents the shape confidence interval, obtained as
$\Ebb[\Sigma](t) \pm 1.96\cdot\StdDev[\Sigma](t)$.
The yellow dots correspond to $\gamma_{\Gamma}(0)$
and $\gamma_{\Sigma_{\mathrm{ref}}(t)}(0)$.
Second row: forward data $u\big(\gamma_{\Sigma}(s),t\big)$.
Third row: inverse data $y_{\text{d}}\big(\gamma_{\Gamma}(s),t\big)$.}
\label{fig:data}
\end{figure}

In Figure~\ref{fig:data}, we summarised the input data for the forward and
inverse problem and the geometry of the pericardium, superimposed with
the CMR images.
To model the uncertainty, we selected the resulting reference shape of the
pericardium $\Sigma_{\text{ref}}(t)$, for $t\in[0,T)$, as the space-time mean
of the random deformation field $\boldsymbol{\chi}$. The covariance of
$\boldsymbol{\chi}$ evaluated at the points $\hat{\bf z}=(\hat{\bf x},t)$ with
$\hat{\bf x}\in\Sigma_{\text{ref}}(t)$ and $\hat{\bf z}^\prime=(\hat{\bf
x}^\prime,t^\prime)$ with $\hat{\bf x}^\prime\in\Sigma_{\text{ref}}(t^\prime)$
was the product of a matrix-valued kernel in space and a scalar kernel in time.
For the spatial component of the covariance, we considered the Mat\'ern kernel
\[
k_{\nu}^{\text{M}}(d)=\sigma^2\frac{2^{1-\nu}}{\Gamma(\nu)}
\left(\sqrt{2\nu}\frac{d}{\rho}\right)^{\nu}
K_{\nu}\left(\sqrt{2\nu}\frac{d}{\rho}\right),
\]
with parameters $\rho=50$, $\sigma^2=4/3$ and different values of the
smoothness index \(\nu\). Note that in the definition of \(k_{\nu}^{\text{M}}\),
the function $\Gamma(\cdot)$ is the gamma function and $K_{\nu}(\cdot)$ is the
modified Bessel function of the second kind. Since we are interested in
modeling the uncertainty due to the noise and the limited resolution of CMR,
the distance $d=d(\xmat,\xmat')$ used in the
definition of the kernel, was measured with the Euclidean norm in $\mathbb{R}^2$.
In contrast, for the distance in
time between $t$ and $t^\prime$, we take into account the periodicity of the
motion of the pericardium. To do so, we first scale the interval $[0,T)$ to
$[0,2\pi)$ with the mapping
\begin{align*}
\eta\colon[0,T)&\rightarrow[0,2\pi),\quad
t\mapsto\frac{2\pi t}{T},
\end{align*}
and then we compute the geodesic distance
\[
\theta(t,t^\prime)=\arccos\Big(\cos\big(\eta(t)-\eta(t^\prime)\big)\Big).
\]
Finally, to model an additional correlation between neighbouring time
slices, we employ the sine power kernel
\[
k^{\text{SP}}(\theta)=1-\big(\sin(\theta/2)\big)^2=\big(\cos(\theta/2)\big)^2=\frac{1}{2}\big(\cos(\theta)+1\big), 
\]
see \cite{Gneiting}.
Employing a tensor product construction, we end up with the
 covariance function
\[
\Cov[\boldsymbol{\chi}](\zmat,\zmat^\prime)=k^{\text{SP}}\big(\theta(t,t^\prime)\big)
\begin{bmatrix}
k_{5/2}^{\text{M}}(\|\xmat(t)-\xmat^\prime(t')\|_2) & 0\\
0 & k_{\infty}^{\text{M}}(\|\xmat(t)-\xmat^\prime(t')\|_2)
\end{bmatrix}.
\]
Note that we assume here that there is no correlation between the
two spatial components of the deformation field.  By construction,
the joint covariance kernel is positive-definite on the space
$\mathbb{R}^2\times \mathbb{S}^1$, where $\mathbb{S}^1$ is the unit circle.

In the forward problem, the pericardial potential $u(\xbf,t)$ was defined
analytically as a $T$-periodic function in the variable $t$.
For convenience, we set $u(\xbf,t) = u\big(\gamma_{\Sigma(t,\omega)}(s),t\big)$,
being $s\in[0,1)$ the normalised curvilinear coordinate.
We simulated a left bundle branch block, that is the extracellular potential
consisted in a propagation from free wall of the right ventricle,
at $s=0$, towards the free wall of the left ventricle, at $s=0.5$.
The propagation took \SI{150}{\ms}, consistent with a long QRS complex.
The specific analytical form of $u(\xbf,t)$ is given in
Appendix~\ref{sec:forwardData}.

To generate the input data for the inverse problem, we computed the solution of
the forward problem on the space-time reference geometry and eventually added
Gaussian noise with zero mean and variance $10^{-8}$, corresponding to
a signal-to-noise ratio of approximately \SI{46}{\decibel}.

Concerning the space discretisation, we considered $n_{\Sigma}=n_{\Gamma}=500$
collocation points. This resulted in a truncated Karhunen-Lo\`eve expansion 
with $K=648$ terms, obtained with a tolerance of $10^{-4}$.  That is, the
parametric dimension of the UQ problem was $648$, thus high-dimensional.
We remark that, since the boundaries are represented by trigonometric
polynomials and all data are smooth functions, the collocation method
converges exponentially.  Resulting in spatial approximation errors for the
chosen number of boundary points, which are already of the order of
the machine precision.

\begin{figure}[ht]
\centering
\begin{tikzpicture}
\begin{loglogaxis}[width=0.5\textwidth,grid=both, ymin= 3e-5, ymax = 3e-3, xmin = 1e2, xmax =2e4,
    legend style={legend pos=north east,font=\small}, ylabel={\small Error}, xlabel ={\small QMC samples $N$}]
\addplot[line width=0.7pt,color=red,mark=o] table[x=N,y=M1]{./resultsRevision/Forward.txt};\addlegendentry{$\mathcal{M}_1$};
\addplot[line width=0.7pt,color=blue,mark=triangle] table[x=N,y=M2]{./resultsRevision/Forward.txt};\addlegendentry{$\mathcal{M}_2$};
\addplot[line width=0.7pt, color=black,dashed] table[x=N, y={create col/linear regression={y=M1}}]{./resultsRevision/Forward.txt};
\xdef\slope{\pgfplotstableregressiona}
\addlegendentry{$N^{\pgfmathprintnumber{\slope}}$};
\addplot[line width=0.7pt, color=black,dashed,] table[x=N, y={create col/linear regression={y=M2}}]{./resultsRevision/Forward.txt};
\xdef\slopep{\pgfplotstableregressiona}
\addlegendentry{$N^{\pgfmathprintnumber{\slopep}}$};
\end{loglogaxis}
\end{tikzpicture}
\caption{Convergence plot of the first ($\mathcal{M}_1$) and second ($\mathcal{M}_2$) moment of the forward solution.}
\label{fig:Forward}
\end{figure}
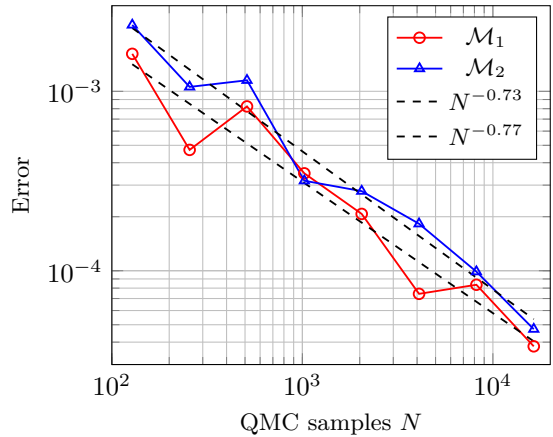

Since we wish to employ the sparse quadrature to estimate our quantities of
interest, we first numerically tested its convergence by comparing it to the
quasi-Monte Carlo method based on the Halton set, see~\cite{Caf98}.  We tested
the convergence of the first and second moment of the forward and inverse
solution at $t=\SI{189}{\ms}$.  The time-independent problem yielded
a truncated Karhunen-Lo\`eve expansion of $K=101$ terms. To validate the
applicability of the sparse quadrature, we consider the approximation
error of the moments by a comparison to the quasi-Monte Carlo method
based on Halton points. The accuracy of approximately \(5\cdot 10^-5\) was
achieved by using \num{17799} quadrature points within the sparse quadrature.
Figure~\ref{fig:Forward} shows the convergence plot of the forward solution.
We experimentally observed a convergence rate of 0.75. The obtained result
indeed corroborates that the forward problem is smooth
in the parametric space.

For the inverse problem we tested the regularisations reported in
Section~\ref{sec:inverse}.  The regularisation parameter was determined with
the L-curve method, see \cite{mueller2012linear}, on the space-time reference
geometry.  For each regularisation, we then chose the maximal regularisation
parameter over all the time steps.  This might have led to over-regularisation
for some time steps, but it has the benefit of limiting the oscillations in the
numerical solution of the inverse problem.  The values of the chosen
regularisation parameters $\lambda$ are reported in Table~\ref{table:param}.

\begin{table}
\begin{center}
\begin{tabular}{ c | c | c | c | c }
 & Zero order Tikhonov & First order Tikhonov & $H^{1/2}$ & Total Variation \\
 \hline
 $\lambda$ & $10^{-6}$ & $10^{-3}$ & $10^{-5}$ & $10^{-5}$ 
\end{tabular}
\end{center}
\caption{regularisation parameter choice for the zero order Tikhonov,
first order Tikhonov, $H^{1/2}$ and Total Variation regularisations.}
\label{table:param}
\end{table}

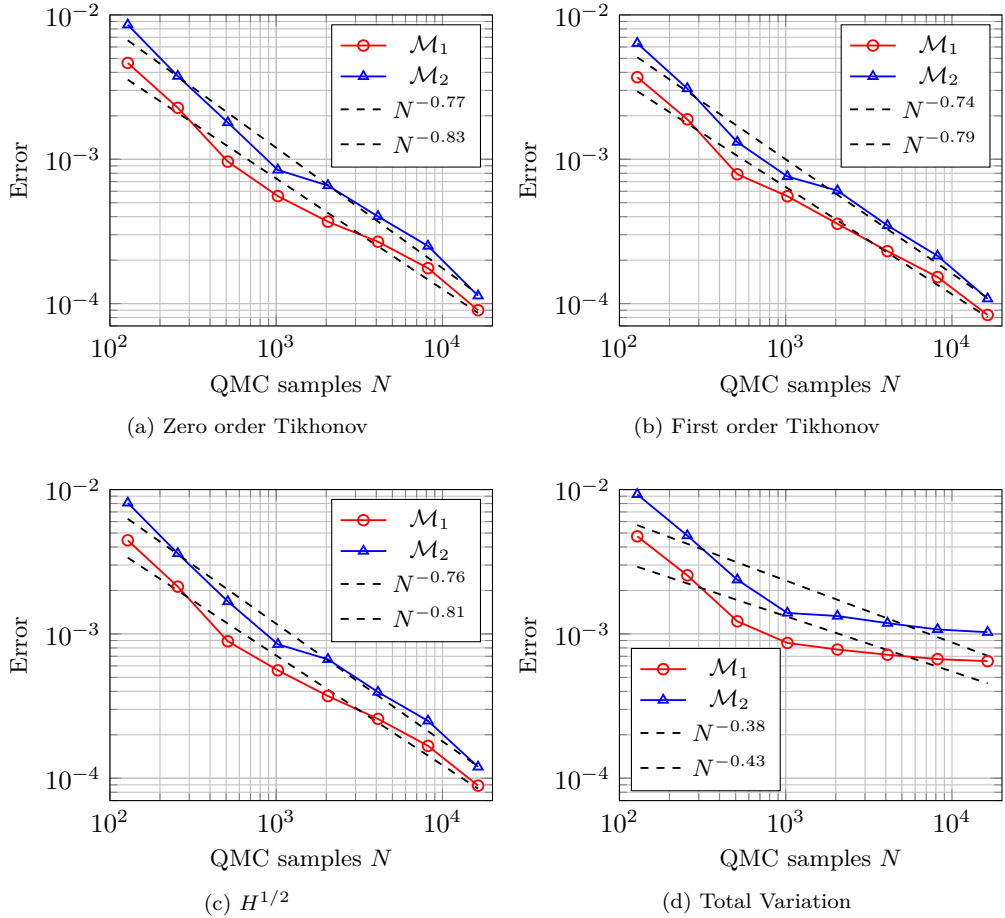
\begin{figure}[ht]
\subfloat[][Zero order Tikhonov]{
\begin{tikzpicture}
\begin{loglogaxis}[width=0.45\textwidth,grid=both, ymin= 7e-5, ymax = 1e-2, xmin = 1e2, xmax =2e4,
    legend style={legend pos=north east,font=\small}, ylabel={\small Error}, xlabel ={\small QMC samples $N$}]
\addplot[line width=0.7pt,color=red,mark=o] table[x=N,y=M1]{./resultsRevision/InvZeroTik.txt};\addlegendentry{$\mathcal{M}_1$};
\addplot[line width=0.7pt,color=blue,mark=triangle] table[x=N,y=M2]{./resultsRevision/InvZeroTik.txt};\addlegendentry{$\mathcal{M}_2$};
\addplot[line width=0.7pt, color=black,dashed] table[x=N, y={create col/linear regression={y=M1}}]{./resultsRevision/InvZeroTik.txt};
\xdef\slope{\pgfplotstableregressiona}
\addlegendentry{$N^{\pgfmathprintnumber{\slope}}$};
\addplot[line width=0.7pt, color=black,dashed,] table[x=N, y={create col/linear regression={y=M2}}]{./resultsRevision/InvZeroTik.txt};
\xdef\slopep{\pgfplotstableregressiona}
\addlegendentry{$N^{\pgfmathprintnumber{\slopep}}$};
\end{loglogaxis}
\end{tikzpicture}
}
\subfloat[][First order Tikhonov]{
\begin{tikzpicture}
\begin{loglogaxis}[width=0.45\textwidth,grid=both, ymin= 7e-5, ymax = 1e-2, xmin = 1e2, xmax =2e4,
    legend style={legend pos=north east,font=\small}, ylabel={\small Error}, xlabel ={\small QMC samples $N$}]
\addplot[line width=0.7pt,color=red,mark=o] table[x=N,y=M1]{./resultsRevision/InvFirstTik.txt};\addlegendentry{$\mathcal{M}_1$};
\addplot[line width=0.7pt,color=blue,mark=triangle] table[x=N,y=M2]{./resultsRevision/InvFirstTik.txt};\addlegendentry{$\mathcal{M}_2$};
\addplot[line width=0.7pt, color=black,dashed] table[x=N, y={create col/linear regression={y=M1}}]{./resultsRevision/InvFirstTik.txt};
\xdef\slope{\pgfplotstableregressiona}
\addlegendentry{$N^{\pgfmathprintnumber{\slope}}$};
\addplot[line width=0.7pt, color=black,dashed,] table[x=N, y={create col/linear regression={y=M2}}]{./resultsRevision/InvFirstTik.txt};
\xdef\slopep{\pgfplotstableregressiona}
\addlegendentry{$N^{\pgfmathprintnumber{\slopep}}$};
\end{loglogaxis}
\end{tikzpicture}
}\\
\subfloat[][$H^{1/2}$]{
\begin{tikzpicture}
\begin{loglogaxis}[width=0.45\textwidth,grid=both, ymin= 7e-5, ymax = 1e-2, xmin = 1e2, xmax =2e4,
    legend style={legend pos=north east,font=\small}, ylabel={\small Error}, xlabel ={\small QMC samples $N$}]
\addplot[line width=0.7pt,color=red,mark=o] table[x=N,y=M1]{./resultsRevision/InvHhalf.txt};\addlegendentry{$\mathcal{M}_1$};
\addplot[line width=0.7pt,color=blue,mark=triangle] table[x=N,y=M2]{./resultsRevision/InvHhalf.txt};\addlegendentry{$\mathcal{M}_2$};
\addplot[line width=0.7pt, color=black,dashed] table[x=N, y={create col/linear regression={y=M1}}]{./resultsRevision/InvHhalf.txt};
\xdef\slope{\pgfplotstableregressiona}
\addlegendentry{$N^{\pgfmathprintnumber{\slope}}$};
\addplot[line width=0.7pt, color=black,dashed,] table[x=N, y={create col/linear regression={y=M2}}]{./resultsRevision/InvHhalf.txt};
\xdef\slopep{\pgfplotstableregressiona}
\addlegendentry{$N^{\pgfmathprintnumber{\slopep}}$};
\end{loglogaxis}
\end{tikzpicture}
}
\subfloat[][Total Variation]{
\begin{tikzpicture}
\begin{loglogaxis}[width=0.45\textwidth,grid=both, ymin= 7e-5, ymax = 1e-2, xmin = 1e2, xmax =2e4,
    legend style={legend pos=south west,font=\small}, ylabel={\small Error}, xlabel ={\small QMC samples $N$}]
\addplot[line width=0.7pt,color=red,mark=o] table[x=N,y=M1]{./resultsRevision/InvTotVar.txt};\addlegendentry{$\mathcal{M}_1$};
\addplot[line width=0.7pt,color=blue,mark=triangle] table[x=N,y=M2]{./resultsRevision/InvTotVar.txt};\addlegendentry{$\mathcal{M}_2$};
\addplot[line width=0.7pt, color=black,dashed] table[x=N, y={create col/linear regression={y=M1}}]{./resultsRevision/InvTotVar.txt};
\xdef\slope{\pgfplotstableregressiona}
\addlegendentry{$N^{\pgfmathprintnumber{\slope}}$};
\addplot[line width=0.7pt, color=black,dashed,] table[x=N, y={create col/linear regression={y=M2}}]{./resultsRevision/InvTotVar.txt};
\xdef\slopep{\pgfplotstableregressiona}
\addlegendentry{$N^{\pgfmathprintnumber{\slopep}}$};
\end{loglogaxis}
\end{tikzpicture}
}
\caption{Convergence plot of the first ($\mathcal{M}_1$) and second ($\mathcal{M}_2$) moment of the inverse solution with (a) zero order Tikhonov, (b) first order Tikhonov, (c) $H^{1/2}$ and (d) Total Variation regularisation.}
\label{fig:Inverse}
\end{figure}

Figure~\ref{fig:Inverse} shows the convergence plots of the inverse solutions.
We can observe convergence towards the reference moments for the zero order
Tikhonov, the first order Tikhonov and the $H^{1/2}$ regularisations.  Instead,
the curves resulting from the Total Variation regularisation flatten out,
meaning that there is no convergence and suggesting that this regularisation is
not sufficiently smooth with respect to the random parameter, and therefore
inadequate for the sparse quadrature approach.
In contrast, the other regularisations are to be suitable for the sparse
quadrature approach.

\begin{figure}
\includegraphics[width=\textwidth, keepaspectratio]{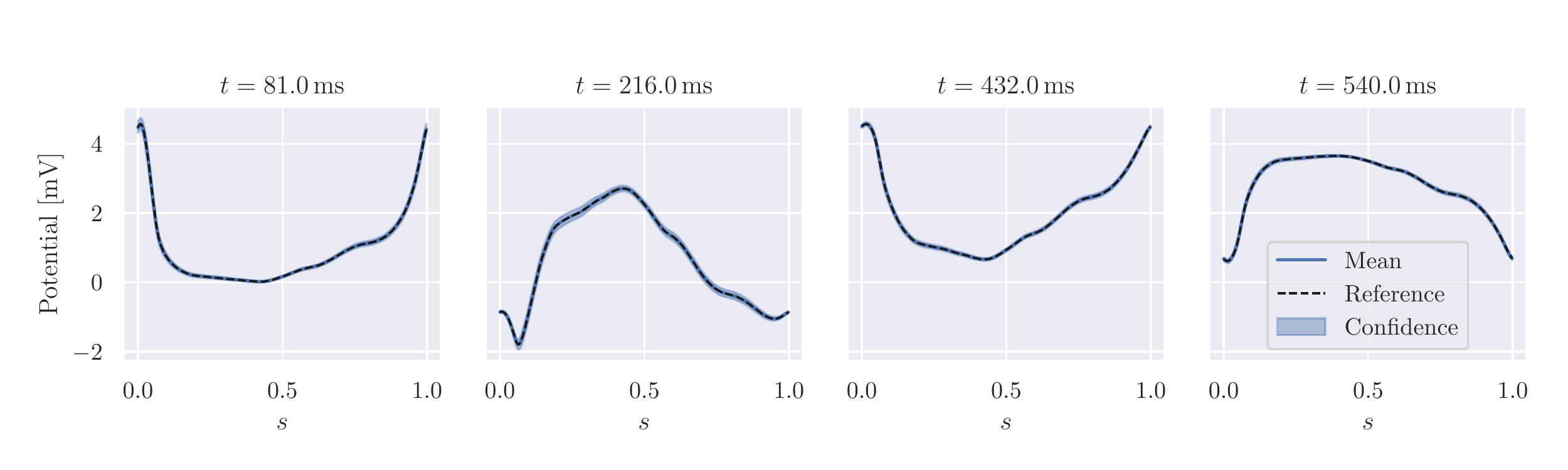}
\includegraphics[width=\textwidth, keepaspectratio]{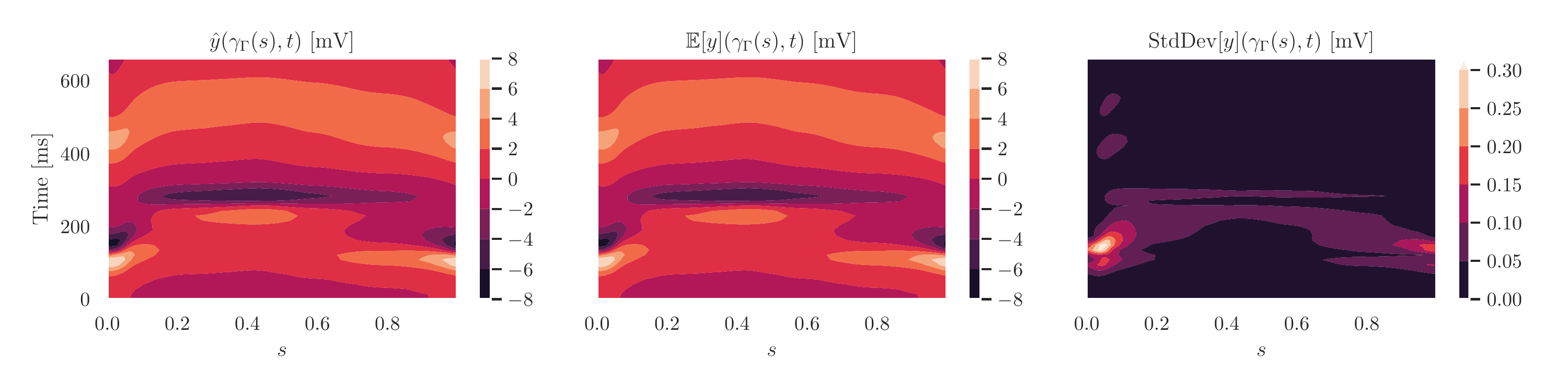}
\caption{Solution of the forward problem. First row: expectation and confidence
interval (blue), and solution with reference geometry (dashed black).
Second row: contour plots in space-time of the reference solution, expected
solution and standard deviation.}
\label{fig:resultsForward}
\end{figure}

Next, we present the corresponding numerical results for the inverse problem,
when employing the $H^{1/2}$ regularisation. We estimate the mean
and the standard deviation of the forward and inverse solution.
The computed quantities of interest for the forward problem are shown in Figure~\ref{fig:resultsForward}.
We can observe that the expectation was very close to the chest potential
computed on the space-time reference geometry.  Moreover, the standard
deviation is small and it mainly affects the regions with higher amplitude
during the depolarisation phase.  In particular, the standard deviation is
higher close to $\gamma_{\Gamma}(0)$. The reason for this phenomenon is that
points in the vicinity of $\gamma_{\Gamma}(0)$ are close to the heart, see Figure~\ref{fig:data}. 
Therefore the effect of the shape uncertainty on the forward solution
is limited.

\begin{figure}
\includegraphics[width=\textwidth, keepaspectratio]{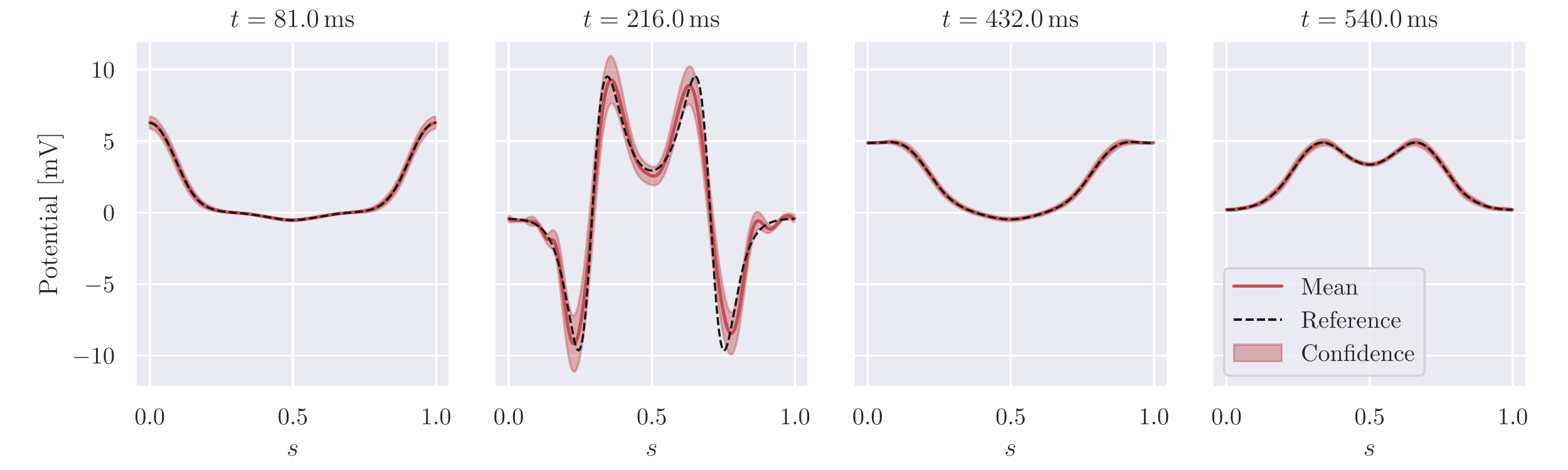}
\includegraphics[width=\textwidth, keepaspectratio]{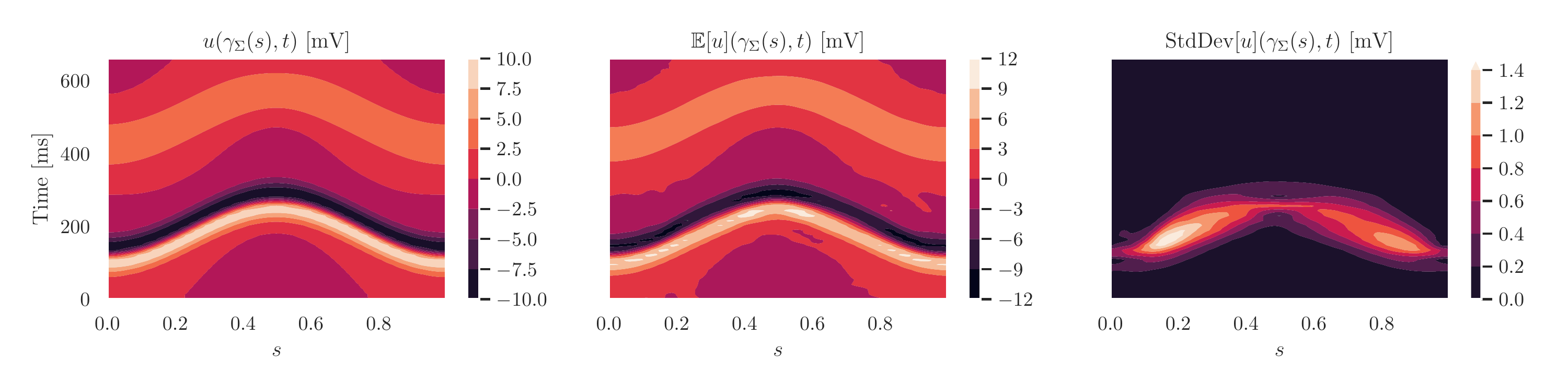}
\caption{Solution of the inverse problem. First row: expectation and confidence
interval (red), and reference solution (dashed black).
Second row: contour plots in space-time of the reference solution, expected
solution and standard deviation.}
\label{fig:resultsInverse}
\end{figure}

The computed quantities of interest for the inverse problem are shown in Figure~\ref{fig:resultsInverse}.
We observe that the expectation deviates slightly from the pericardial
potential and shows some oscillations in the depolarisation phase. 
This is due to the ill-posedness of the inverse problem and the effect of the
regularisation. Indeed, the Tikhonov regularisation, especially at
lower order, tends to introduce oscillations in the solution, if the regularisation
parameter is not sufficiently large.  On the other hand, a strong regularisation
yields smoothed out gradients, thus less accurate inverse solution.
The $H^{1/2}$ regularisation might therefore
introduce oscillations as well.
In the bottom row of Figure~\ref{fig:resultsInverse}, the standard deviation is
larger than in the forward problem and the shape uncertainty mainly affects
regions with large gradients.  In the depolarisation phase, the potential is
steeper and, consequently, the standard deviation is larger.
This fact may have a consequence on the quality of derived quantities such as
the activation map, usually obtained from the point with largest negative
deflection in the signals.
In conclusion, the effect of the shape uncertainty is much more relevant for
the inverse problem than for the forward problem.

\section{Conclusions}

In this work we have considered the forward and inverse problem of
electrocardiography in the presence of space-time shape uncertainties. The
high parametric dimensionality of the problem, along with the non-linear map
between the input and output uncertainties, renders the problem challenging
from the mathematical and the numerical perspective. 

To address space-time shape uncertainties, we have suggested a model by random
space-time deformation fields with a periodic-in-time covariance kernel. This
approach resulted in a high-dimensional uncertainty quantification problem. To
make this dimensionality feasible for numerical computations, we rely on the
one hand on a low-rank representation of the random deformation field and
efficient quadrature techniques on the other hand. The resulting spatial
problem for each quadrature point in the parameter was addressed by a boundary
integral formulation in combination with a highly accurate and fast collocation
method. The numerical results indicate that shape uncertainties affect the
solution much more strongly in the inverse problem than in the forward one,
as a consequence of the ill-posedness.  Moreover, it was observed that the
regularisation of the inverse solution may affect the regularity of the
solution with respect to the stochastic parameter, with possible limitations in
the convergence order of the quadrature method. Especially, the total variation
regularisation showed poor performance in this respect, while the $H^{1/2}$
regularisation is suitable for this problem, both in terms of regularity and
quality of the reconstruction.

Future research directions include the transition to 3-dimensional models and the
generalisation to piecewise-constant conductivities in the torso. In this case,
it is still possible to leverage on the boundary integral formulation. From the
mathematical perspective, a general result on the parametric regularity of
the inverse problem is still missing. Such a result would help to shed some light
on the sub-optimal convergence rate observed in our experiments and the apparent
lack of parametric regularity of the inverse problem with TV regularisation.
Finally, from the clinical perspective, we shall further investigate how
segmentation uncertainty may be modelled.  As a matter of fact, the uncertainty
likely results from multiple sources. Besides noise, CMR images are affected
by breath holding, which may limit the volume of the heart and spatio-temporal
mis-registration. The segmentation process may amplify the acquisition
uncertainty, depending on the method.  Importantly, the segmentation is
sometimes a non-smooth operation, in the sense that small perturbation in the
images could yield large and possibly topological changes in the contour.
Therefore, special care will be needed in order to correctly model the
uncertainty from clinical data.

\paragraph{Acknowledgements}
This work was financially supported by the Theo Rossi di Montelera Foundation,
the Metis Foundation Sergio Mantegazza,
the Fidinam Foundation,
the Horten Foundation and the CSCS—Swiss National Supercomputing Centre
production grant s1074.

\appendix
\section{Collocation method}\label{sec:Nystrom}
Considering the parameterisations $\gamma_{\Gamma}\colon [0,1)\rightarrow\Gamma$
and $\gamma_{\Sigma}\colon[0,1)\rightarrow\Sigma$, the Dirichlet-to-Neumann
map \eqref{eq:BIE} in dimension $d=2$ reads, for $\Phi\in\{\Gamma,\Sigma\}$
and $s\in[0,1)$,
\begin{align*}
-\frac{1}{4\pi}\sum_{\Psi\in\{\Gamma,\Sigma\}}\int_0^1
k^{\mathcal{V}}_{\Phi,\Psi}(s,r)\tilde{\rho}_{1,\Psi}(r)\dd
r=\frac{1}{2}\rho_{0,\Phi}(s)+\frac{1}{2\pi}\sum_{\Psi\in\{\Gamma,\Sigma\}}\int_0^1k^{\mathcal{K}}_{\Phi,\Psi}(s,r) \rho_{0,\Psi}(r)\dd r,
\end{align*}
where
\begin{align*}
k^{\mathcal{V}}_{\Phi,\Psi}(s,r)\isdef{\log\big\|\gamma_{\Phi}(s)
-\gamma_{\Psi}(r)\big\|_2^2},\quad
k^{\mathcal{K}}_{\Phi,\Psi}(s,r)\isdef\frac{\big\langle\gamma_{\Phi}(s)
-\gamma_{\Psi}(r),\nbf_{r,\Psi}\big\rangle_2}{\big\|\gamma_{\Phi}(s)
- \gamma_{\Psi}(r)\big\|_2^2}\big\|\gamma_{\Psi}^\prime(r)\big\|_2
\end{align*}
and
\[
\rho_{0,\Phi}(s)\isdef y\big(\gamma_{\Phi}(s)\big),\quad
\tilde{\rho}_{1,\Phi}(s)\isdef \frac{\partial y}{\partial
\nbf_{s,\Phi}}\big(\gamma_{\Phi}(s)\big)\big\|\gamma_{\Phi}^\prime(s)\big\|_2.
\]
We consider the $n_{\Gamma}$-points discretisation of $\Gamma$ and the
$n_{\Sigma}$-points discretisation of $\Sigma$, where $n_{\Gamma}$
and $n_{\Sigma}$ are even numbers, and, for $\Phi\in\{\Gamma,\Sigma\}$,
we represent $\rho_{0,\Phi}$ and $\tilde{\rho}_{1,\Phi}$ as a linear
combination of the trigonometric Lagrange polynomials $L_j$,
for $j=0,\ldots,n_{\Phi}-1$, i.e.,
\[
\rho_{0,\Phi}(s)=\sum_{j=0}^{n_{\Phi}-1}\rho_{0,j}^{\Phi}L_j(s)\quad\text{and}\quad
\tilde{\rho}_{1,\Phi}(s)=\sum_{j=0}^{n_{\Phi}-1}\tilde{\rho}_{1,j}^{\Phi}L_j(s).
\]
Note that, considering the discretisation $s_i=i/n_{\Phi}$ for $i=0,...,n_{\Phi}-1$, it holds $L_j(s_i)=\delta_{i,j}$. Using the representation of $\rho_{0,\Phi}$ and $\tilde{\rho}_{1,\Phi}$, we end up with
\begin{align*}
-\frac{1}{4\pi}\sum_{\Psi\in\{\Gamma,\Sigma\}}\sum_{j=0}^{n_{\Psi}-1}&\left(\int_0^1 k^{\mathcal{V}}_{\Phi,\Psi}(s_i,r)L_j(r)\dd r\right)\tilde{\rho}_{1,j}^{\Psi}\\
&=\frac{1}{2}\rho_{0,i}^{\Phi}+\frac{1}{2\pi}\sum_{\Psi\in\{\Gamma,\Sigma\}}\sum_{j=0}^{n_{\Psi}-1}\left(\int_0^1k^{\mathcal{K}}_{\Phi,\Psi}(s_i,r)L_j(r)\dd r\right)\rho_{0,j}^{\Psi}
\end{align*}
for $i=0,\ldots,n_{\Phi}-1$ and $\Phi\in\{\Gamma,\Sigma\}$.
The goal is now to discretise the above equation. To do so,
we consider the discretisation $r_j=j/n_{\Psi}$ for $j=0,\ldots,n_{\Psi}-1$
and we distinguish between four different kinds of matrices.

In the case $\Psi\neq\Phi$, the kernel functions do not exhibit singularities.
Hence, we may directly employ the trapezoidal rule
to obtain the collocation approximation of the corresponding boundary
integral operators.
Using the trapezoidal rule we have
\begin{align*}
-\frac{1}{4\pi}\sum_{j=0}^{n_{\Psi}-1}&
\left(\int_0^1 k^{\mathcal{V}}_{\Phi,\Psi}(s_i,r)L_j(r)\dd r\right)
\tilde{\rho}_{1,j}^{\Psi}\approx-\frac{1}{4\pi n_{\Psi}}
\sum_{j=0}^{n_{\Psi}-1}k^{\mathcal{V}}_{\Phi,\Psi}(s_i,r_j)
\tilde{\rho}_{1,j}^{\Psi},
\end{align*}
therefore we define the matrix $\Vbf_{\Phi\Psi}
\in\mathbb{R}^{n_{\Phi}\times n_{\Psi}}$
whose entries, for $i=0,\ldots,n_{\Phi}-1$ and
$j=0,\ldots,n_{\Psi}-1$, are
\[
(\Vbf_{\Phi\Psi})_{i,j}\isdef-\frac{1}{4\pi n_{\Psi}}
k^{\mathcal{V}}_{\Phi,\Psi}(s_i,r_j).
\]
Similarly, by the trapezoidal rule we have
\begin{align*}
\frac{1}{2\pi}\sum_{j=0}^{n_{\Psi}-1}&
\left(\int_0^1 k^{\mathcal{K}}_{\Phi,\Psi}(s_i,r)L_j(r)\dd r\right)
\rho_{0,j}^{\Psi}\approx\frac{1}{2\pi n_{\Psi}}\sum_{j=0}^{n_{\Psi}-1}
k^{\mathcal{K}}_{\Phi,\Psi}(s_i,r_j)\rho_{0,j}^{\Psi}.
\end{align*}
Hence, we define the matrix
$\Kbf_{\Phi\Psi}\in\mathbb{R}^{n_{\Phi}\times n_{\Psi}}$ whose entries, for
$i=0,\ldots,n_{\Phi}-1$ and $j=0,\ldots,n_{\Psi}-1$, are
\[
(\Kbf_{\Phi\Psi})_{i,j}\isdef\frac{1}{2\pi n_{\Psi}}
k^{\mathcal{K}}_{\Phi,\Psi}(s_i,r_j).
\]

If $\Psi=\Phi$ the kernel functions exhibit a singularity for \(s=r\).
This singularity can be dealt with as follows: 
We split
\[
k^{\mathcal{V}}_{\Phi,\Phi}(s,r)=k^{\mathcal{V},1}_{\Phi,\Phi}(s,r)
+k^{\mathcal{V},2}_{\Phi,\Phi}(s,r),
\]
where
\[
k^{\mathcal{V},1}_{\Phi,\Phi}(s,r)
\isdef\log\frac{\big\|\gamma_{\Phi}(s)-\gamma_{\Phi}(r)
\big\|_2^2}{4\sin^2\big(\pi(s-r)\big)}
\]
and
\[
k^{\mathcal{V},2}_{\Phi,\Phi}(s,r)\isdef\log\Big(4\sin^2\big(\pi(s-r)\big)\Big).
\]
Then, for the first term, using again the trapezoidal rule, we arrive at
\begin{align*}
-\frac{1}{4\pi}\sum_{j=0}^{n_{\Phi}-1}&\left(\int_0^1 k^{\mathcal{V},1}_{\Phi,\Phi}(s_i,r)L_j(r)\dd r\right)\tilde{\rho}_{1,j}^{\Phi}\approx-\frac{1}{4\pi n_{\Phi}}\sum_{j=0}^{n_{\Phi}-1}k^{\mathcal{V},1}_{\Phi,\Phi}(s_i,r_j)\tilde{\rho}_{1,j}^{\Phi}.
\end{align*}
Since
\begin{align*}
\lim_{s\rightarrow r}k^{\mathcal{V},1}_{\Phi,\Phi}(s,r)=\log\left\|\frac{\gamma_{\Phi}^{\prime}(r)}{2\pi}\right\|_2^2,
\end{align*}
we set
\begin{align*}
k^{\mathcal{V},1}_{\Phi,\Phi}(s_i,r_j)\isdef\log\left\|\frac{\gamma_{\Phi}^{\prime}(r_j)}{2\pi}\right\|_2^2
\end{align*}
 for \(s_i= r_j\) and we define the matrix $\Vbf_{\Phi\Phi}^1\in\mathbb{R}^{n_{\Phi}\times n_{\Phi}}$ with
\[
(\Vbf_{\Phi\Phi}^1)_{i,j}\isdef-\frac{1}{4\pi n_{\Phi}}k^{\mathcal{V},1}_{\Phi,\Phi}(s_i,r_j)
\quad\text{for }i,j=0,...,n_{\Phi}-1.
\]
Since $n_{\Phi}$ is an even number, we have $n_{\Phi}=2m_{\Phi}$ and for the second term it holds
\begin{align*}
-\frac{1}{4\pi}\sum_{j=0}^{n_{\Phi}-1}&\left(\int_0^1 k^{\mathcal{V},2}_{\Phi,\Phi}(s_i,r)L_j(r)\dd r\right)\tilde{\rho}_{1,j}^{\Phi}\approx -\frac{1}{4\pi}\sum_{j=0}^{n_{\Phi}-1}R_j(s_i)\tilde{\rho}_{1,j}^{\Phi},
\end{align*}
where
\[
R_j(s_i)=-\frac{1}{m_{\Phi}}\bigg(\sum_{k=1}^{m_{\Phi}-1}\frac{1}{k}\cos\big(2\pi k(s_i-r_j)\big)+\frac{1}{n_{\Phi}}\cos\big(2\pi m_{\Phi}(s_i-r_j)\big)\bigg),
\]
cp. \cite{kress1989linear}.
This gives rise to the matrix $\Vbf_{\Phi\Phi}^2
\in\mathbb{R}^{n_{\Phi}\times n_{\Phi}}$ with
\[
(\Vbf_{\Phi\Phi}^2)_{i,j}\isdef-\frac{1}{4\pi}R_j(s_i)\quad\text{for }
i,j=0,\ldots,n_{\Phi}-1.
\]
Finally, we set
\[
\Vbf_{\Phi\Phi}\isdef \Vbf_{\Phi\Phi}^1
+\Vbf_{\Phi\Phi}^2\in\mathbb{R}^{n_{\Phi}\times n_{\Phi}}.
\]

For the double layer operator, we have
\begin{align*}
\frac{1}{2\pi}\sum_{j=0}^{n_{\Phi}-1}&
\left(\int_0^1 k^{\mathcal{K}}_{\Phi,\Phi}(s_i,r)L_j(r)\dd r\right)
\rho_{0,j}^{\Phi}\approx\frac{1}{2\pi n_{\Phi}}\sum_{j=0}^{n_{\Phi}-1}
k^{\mathcal{K}}_{\Phi,\Phi}(s_i,r_j)\rho_{0,j}^{\Phi}.
\end{align*}
Since
\begin{align*}
\lim_{s\rightarrow r}k^{\mathcal{K}}_{\Phi,\Phi}(s,r)
=\frac{\big\langle\gamma_{\Phi}^{\prime\prime}(r),\nbf_{r,\Phi}\big\rangle_2}{2\big\|
\gamma_{\Phi}^{\prime}(r)\big\|_2},
\end{align*}
we set
\begin{align*}
k^{\mathcal{K}}_{\Phi,\Phi}(s_i,r_j)
\isdef\frac{\big\langle\gamma_{\Phi}^{\prime\prime}(r_j),\nbf_{r_j,\Phi}\big\rangle_2}
{2\big\|\gamma_{\Phi}^{\prime}(r_j)\big\|_2}
\end{align*}
for \(s_i=r_j\) and we define the matrix
$\Kbf_{\Phi\Phi}\in\mathbb{R}^{n_{\Phi}\times n_{\Phi}}$ with
\[
(\Kbf_{\Phi\Phi})_{i,j}\isdef\frac{1}{2\pi n_{\Phi}}
k^{\mathcal{K}}_{\Phi,\Phi}(s_i,r_j)\quad\text{for }
i,j=0,\ldots,n_{\Phi}-1.
\]

\section{Forward problem data}\label{sec:forwardData}
The evolution in time of the pericardial potential at the point
$\gamma_{\Sigma(t,\omega)}(0)$ from which the stimulus is delivered is given,
for $t\in[0,T)$, by
\[
u\big(\gamma_{\Sigma(t,\omega)}(0),t\big)=u_{\text{dep}}(t)+u_{\text{rep}}(t),
\]
with
\[
u_{\text{dep}}(t)=-25\frac{\tanh\Big(2\big((t/T-0.18)-\lfloor 0.5+(t/T-0.18)
\rfloor\big)/0.1\Big)}{\cosh\Big(2\big((t/T-0.18)-\lfloor 0.5+(t/T-0.18)
\rfloor\big)/0.1\Big)^2}
\]
and
\[
u_{\text{rep}}(t)=\frac{25}{2\sqrt{2\pi}}\Big(\exp\big(-100(t/T-0.63)^2\big)
+\exp\big(-100(t/T+0.37)^2\big)\Big).
\]
For the other points on the pericardium $\Sigma(t,\omega)$, the potential is
shifted in time according to the arrival time of the stimulus. For $s\in[0,1)$,
the shift is given by
\[
\delta(s)=0.22\big(\cos(2\pi s-\pi)+1\big)/2
\]
and the pericardial pericardial is given by
\[
u\big(\gamma_{\Sigma(t,\omega)}(s),t\big)=u_{\text{dep}}\big(t-\delta(s)T\big)
+u_{\text{rep}}\big(t-\delta(s)T\big).
\]


\printbibliography

\end{document}